\documentclass[11pt,letterpaper,reqno]{amsart}
\usepackage{amsfonts}
\usepackage{amssymb}
\usepackage{amsmath}

\allowdisplaybreaks

\newtheorem{theorem}{Theorem}[section]
\newtheorem{lemma}[theorem]{Lemma}
\newtheorem{proposition}[theorem]{Proposition}
\newtheorem{corollary}[theorem]{Corollary}

\theoremstyle{definition}

\theoremstyle{plain}

\numberwithin{equation}{section}

\def\rwr{\mathrel{\mathrm{wr}}}
\newcommand{\p}{\varphi} 
\newcommand{\ga}{\gamma}
\newcommand{\Ac}{\mathcal A}
\newcommand{\CC}{\mathcal C}
\newcommand{\GG}{\mathcal G}
\newcommand{\ol}[1]{\overline{#1}}
\newcommand{\olg}[2]{\pi_{#2}(\ol{g}_{#1})}
\newcommand{\olgt}[2]{\olg{#1}{#2}^T}
\newcommand{\og}[1]{\olg{#1}{k}}
\newcommand{\ogt}[1]{\og{#1}^T}
\newcommand{\Z}{\mathbb{Z}}
\newcommand{\N}{\mathbb{N}}
\newcommand{\C}{\mathbb{C}} 
\newcommand{\R}{\mathbb{R}} 
\newcommand{\tsf}[1]{{#1}(\lambda - \mu) I}

\newcommand{\Ga}{\Gamma}
\newcommand{\inv}{^{-1}}
\def\Sch#1#2#3{\ensuremath{\mathsf{Sch}(#1,#2,#3)}}
\def\Sp#1{\ensuremath{\mathsf{Sp}(#1)}}

\begin{document}
    \title{The Spectra of Lamplighter Groups and Cayley Machines}

\author{Mark Kambites}\address{Fachbereich Mathematik/Informatik
\\  Universit\"at Kassel \\  34109 Kassel\\ Germany}\author{Pedro
V.~Silva}\address{Departamento de Matem\'atica Pura\\
Faculdade de Ci\^encias da Universidade do Porto\\ Rua do Campo
Alegre  687 \\ 4169-007 Porto\\ Portugal}
\author{Benjamin Steinberg}\address{School of Mathematics and
  Statistics\\ Carleton
University\\ 1125 Colonel By Drive\\ Ottawa\\ Ontario  K1S 5B6 \\
Canada}\email{kambites@theory.informatik.uni-kassel.de}
\email{pvsilva@fc.up.pt}
\email{bsteinbg@math.carleton.ca}
\thanks{This work was performed while the first author was at Carleton
  University, funded by the Leverhulme Trust. The second and
  third named authors
gratefully acknowledge the support of NSERC and were also
supported in part by the FCT and POCTI approved project
POCTI/MAT/37670/2001 in participation with the European Community
Fund FEDER.  The second author was also supported in part by FCT 
through the \emph{Centro de
  Matem\'atica da Universidade do Porto}}

\begin{abstract}
We calculate the spectra and spectral measures associated to
random walks on restricted wreath products $G\rwr \Z$, with $G$ a
finite group by calculating the Kesten-von Neumann-Serre spectral measures for
the random walks on Schreier graphs of certain groups generated by
automata. This generalises the work of Grigorchuk and \.{Z}uk on
the lamplighter group.  In the process we characterise when the usual
spectral measure for a group generated by automata coincides with the
Kesten-von 
Neumann-Serre spectral measure.  
\end{abstract}
\maketitle
\section{Introduction}

The systematic study of random walks on discrete non-abelian groups began with
the seminal work of Kesten \cite{Kesten} and has since become an
important area of mathematics; it links such diverse fields as
probability theory, group theory, geometry and analysis. Many of the
properties of a random walk are encapsulated in the spectrum and spectral
measure of the associated Markov operators. However, the computation of
these is in general quite hard; there are as yet very few examples of complete
computations of the spectra \cite{Kesten,survey,zuk,GZlamp} and still fewer
of the spectral measure \cite{survey,GZlamp}.

At the same time, there has been increasing interest in the class of
\textit{automata groups}, or groups generated by transformations defined
by finite automata \cite{Gecseg, Horejs}. The study
of such groups has given many insights into
group theory, and led to the solution of a number of long-standing open
problems. The study of random walks on automata groups was initiated by
Bartholdi and Grigorchuk \cite{Hecke}; they introduced methods to calculate
the spectra and spectral measures associated to random walks of automata
groups on the Schreier graphs with respect to their parabolic subgroups.

These techniques were employed by
Grigorchuk and \.{Z}uk \cite{GZlamp}, who computed the spectrum and
spectral measure of the Markov operator associated to a random walk on the
lamplighter group
$$\Z/2\Z\rwr \Z = \left(\bigoplus_{\Z}\Z / 2\Z\right)\rtimes \Z,$$
 by realizing it as the group generated by a
two-state automaton. They showed that the spectral measure was discrete
-- a previously unseen phenomenon -- and subsequently, in collaboration
with Linnell and Schick, used their computation to obtain a negative
answer to a strong form of the Atiyah conjecture concerning $L^2$-Betti
numbers \cite{counterexample}.

Subsequently, Dicks and Schick~\cite{Dicks} computed the spectral
measures for random walks on groups of the form $G\rwr \Z$ with $G$ a
non-trivial finite group using entirely different methods.  There
methods were then generalised by Bartholdi and Woess to random walks
on a more general class of graphs \cite{BartholdiWoess}. 
The main purpose of this paper is to compute, using finite state
automata, the spectral measure for
random walks on wreath product groups of the form $G\rwr \Z$ with $G$
a finite group, thus obtaining a new proof of the results of Dicks and
Schick~\cite{Dicks} along the lines of the case of $G=\Z/2\Z$
considered by Grigorchuk and \.{Z}uk.  This involves showing that the
spectral measure for these groups coincides with the so-called
Kesten-von Neumann-Serre spectral measure, introduced by Grigorchuk
and \.{Z}uk \cite{GZihara}, for 
random walks of certain groups generated by automata on their Schreier graphs
with respect to parabolic subgroups.  In fact, we obtain the spectrum
for a more general class of groups as well as some new computations of
Kesten-von Neumann-Serre spectral measures.  

We also aim to give a coherent, and as self-contained as possible,
treatment of the subject of computing spectral measures of random
walks on automata groups via calculation of their Kesten-von
Neumann-Serre spectral measures.  We give necessary and
sufficient conditions for these measures to coincide, for groups
acting spherically transitively on rooted trees.  So far as we know,
these criteria have not appeared elsewhere.

In \cite{reset} the second and third authors considered a possible connection
between automata groups and Krohn-Rhodes theory \cite{KR,Arbib}, by
studying the automata group $\Ga$ generated by the Cayley machines of
a finite group $G$. In the case that $G$ is abelian, it transpires that
$\Ga$ is actually the wreath product $G\rwr \Z$. Indeed, when $G = \Z/2\Z$
the automaton obtained is exactly that used by Grigorchuk and \.{Z}uk
to generate the lamplighter group. In the non-abelian case, $\Ga$ is
still a semidirect product of a locally finite group with $\Z$ \cite{reset};
however, we shall see below that $\Ga$ is not in general isomorphic to
$G\rwr \Z$.  The walks that Dicks and Schick considered on $G\rwr \Z$
for $G$ non-abelian can be obtained from the abelian case.

In this paper, we compute the spectrum of the Markov operator associated
to a simple random walk on any automata group generated by a Cayley machine,
proving
that it is always the entire interval $[-1,1]$.  We compute also the
so-called Kesten-von Neumann-Serre spectral measure \cite{GZihara} for the
random walks on the Schreier graphs with respect to parabolic subgroups.
Interestingly, this measure turns out always to be discrete, with support
at the points $\lambda_{p,q}=\cos \frac{p}{q}\pi$ where $1\leq p<q$, $(p,q)=1$;
the weight at $\lambda_{p,q}$ is $\frac{(n-1)^2}{n^q-1}$. In the case that
$G$ is abelian, we show that this coincides with the spectral measure of
the corresponding simple random walk.  
We also calculate the Ihara zeta functions \cite{GZihara} of the
Schreier graphs with respect to parabolic subgroups and hence for
the Cayley graph of $G\rwr \Z$, in the case $G$ is abelian.



\section{Markov and Hecke operators, spectral measures,
random walks and Ihara zeta functions} 

In this section, we introduce the fundamental notions from analysis which
motivate the work in this area. These include Markov and Hecke operators,
spectral measures and Ihara zeta functions. We also see how these concepts
relate to the study of random walks on graphs.

\subsection{Markov and Hecke-type operators}
We consider only real Hilbert spaces.  Let $X=(V,E)$ be a
$k$-regular (undirected) graph with vertex set $V$ and edge set
$E$.   We allow multiple edges and loops.  The primary example for
us is where $\Ga$ is a finitely generated group with symmetric
generating set $S$ of size $k$, $P\leq \Ga$ is a subgroup and $X =
\Sch \Ga P S$ is the associated (left) \emph{Schreier} (or coset)
graph. The vertices are the left cosets $\Ga/P$.  The edge set is
$S\times \Ga/P$.  The edge $(s,gP)$ goes from $gP$ to $sgP$.  In
particular, the (left) \emph{Cayley graph} of $\Ga$ is $\Sch \Ga 1
S$.

For $v\in V$, let $E(v)$ be the set of edges incident with $v$.
For each edge $e\in E(v)$, let $o_v(e)$ denote the vertex at the
other end of $e$ from $v$; for loops $o_v(e)$ is taken to be $v$.
The \emph{random walk} or \emph{Markov operator}
\cite{Kesten,GZlamp} on $X$ is the operator $M:\ell^2(V)\to
\ell^2(V)$ given by
\begin{equation}\label{markov}
Mf(v) = \frac{1}{k}\sum_{e\in E(v)} f(o_v(e))
\end{equation}
Here $\ell^2(V)$ is the space of square summable functions from
$V$ to $\R$.  If $\delta_v : V \to \R$ is the characteristic
function of $\lbrace v \rbrace$ for each $v \in V$, and we write
the ``matrix" for $M$ with respect to the basis $\lbrace \delta_v
\rbrace_{v \in V}$, then we obtain the normalized incidence matrix
for $X$.  That is, the matrix coefficient $\langle
M\delta_{v_1},\delta_{v_2}\rangle$ is the probability that an edge
incident on $v_1$ is also incident on $v_2$.   
For the case of a random walk on a Schreier graph \Sch \Ga P S,
one has that $M:\ell^2(\Ga/P)\to \ell^2(\Ga/P)$ is given by
$$Mf(gP) = \frac{1}{|S|}\sum_{s\in S} f(sgP).$$  This is a special
case of a Hecke-type operator \cite{Hecke,GZlamp}.

Let $\Ga$ be a (discrete) group with finite symmetric generating
set $S$ and $\pi:\Ga\to \mathcal B(\mathfrak H)$ be a unitary
representation of $\Ga$ on a Hilbert space $\mathfrak H$; here
$\mathcal B(\mathfrak H)$ denotes the algebra of bounded linear
operators on $\mathfrak H$. The associated Hecke-type operator is
then $H_{\pi}:\mathfrak H\to \mathfrak H$ given by
$$H_{\pi} = \frac{1}{|S|}\sum_{s\in S}\pi (s).$$  Then $H_{\pi}$ is
a self-adjoint operator and $\|H_{\pi}\|\leq 1$ by the triangle
inequality.

Suppose now that $P\leq \Ga$ is a subgroup and let
$\lambda_{\Ga/P}:\Ga\to \mathcal B(\ell^2(\Ga/P))$ be the (left)
quasi-regular representation; so $\lambda_{\Ga/P}(g)f(hP) =
f(g\inv hP)$.  Then $H_{\lambda_{\Ga/P}}$ is precisely the Markov
operator associated to the random walk on \Sch \Ga P S.  The case
$P=1$ is called the \emph{left regular representation}, denoted
$\lambda_\Ga$.

Another important example is the following. Suppose that $\Ga$
acts on a measure space $(X,\mu)$ by measure-preserving
transformations. Then there is an associated unitary
representation $\pi:\Ga\to \mathcal B(L^2(X,\mu))$ (where
$L^2(X,\mu)$ is the space of square-integrable functions on $X$)
given by $\pi(g)f(x) = f(g^{-1}x)$.  The case of interest to us
arises from the action of an automata group on the boundary of a
rooted tree, viewed as a measure space with the product
(Bernoulli) measure.

\subsection{Kesten spectral measures}
The \emph{spectrum \Sp T} of a bounded operator $T$ on
a Hilbert space consists of all real numbers $\lambda$ such that
$T-\lambda$ is not invertible.  Then $\Sp T$ is a closed subset of
the interval $[-\|T\|,\|T\|]$.

We return now to the Markov operator $M$ for a random walk on a
$k$-regular graph $X$.  As in the case of Schreier graphs, one can
verify that $\|M\|\leq 1$.
 Also $M$ is self-adjoint, so it has a
spectral decomposition
\begin{equation}\label{specdecomp}
M= \int_{-1}^1 \lambda \mathrm dE(\lambda)
\end{equation}
where $E$ is the spectral measure \cite{Pedersen}.  That is, $E$
is a projection-valued measure defined on the Borel subsets of
$[-1,1]$, taking values in the projections of $\mathcal
B(\mathfrak H)$.

For those unfamiliar with these notions, if $X$ is a finite graph,
then
\begin{equation}\label{finitedecomp}
E(B) = \sum_{B\cap \Sp M} E_{\lambda}
\end{equation}
where $E_{\lambda}$ is the projection to the eigenspace associated
with $\lambda$ and \eqref{specdecomp} is the usual orthogonal
decomposition of a symmetric matrix into its eigenspaces.

The matrix $\mu^X$ of measures associated with $E$ is given by
\begin{equation}\label{definespectmeasures}
\mu^X_{v_1,v_2}(B) = \langle E(B)\delta_{v_1},\delta_{v_2}\rangle.
\end{equation}
Of particular interest are the diagonal entries $\mu^X_v
:=\mu^X_{v,v}$.  These are called the \emph{Kesten} spectral
measures associated to the random walk \cite{GZihara}.  To explain
the significance of these measures, we remind the reader about the
moments of a measure.  If $\mu$ is a (Borel) probability measure on
$\R$ (we do not necessarily assume the support of the measure is the
whole real line)
and $f:\R\to \mathbb{R}$ is a measurable function, then the
\emph{expected value} of $f$, denoted $E[f]$, (or $E_{\mu}[f]$ if
we want to emphasize the measure) is given by
\begin{equation*}
E[f] = \int_\R f(x)\mathrm d\mu(x).
\end{equation*}  The \emph{$m^{th}$ moment} of $\mu$, denoted
$\mu^{(m)}$ is $E[x^m]$.  The (formal) \emph{moment generating function} is
then the power series
$$\sum_{m=0}^{\infty}\frac{1}{m!}\mu^{(m)}t^m\in
\mathbb{R}[\!{[}t]\!{]}.$$  This is the Taylor expansion about $0$
of the function $M(t) = E[e^{tx}]$ (where the integral is taken
over $x$). Returning to the situation of a simple random walk on a
graph $X$ and the spectral measure $\mu^X_v$, denote by $p_m(v)$
the probability of return to $v$ on the $m^{th}$ step of the
random walk; then $p_m(v)$ is the $m^{th}$ moment of $\mu^X_v$
\cite{Kesten,GZlamp}. In the case of a Schreier graph \Sch \Ga P S
we shall use the notation $\mu^{\Ga/P}$ to denote $\mu^{\Ga/P}_P$.

If the automorphism group of $X$ acts transitively on the vertices
(for example this occurs for \Sch \Ga P S when $P$ is normal) then
the Kesten measures are independent of the chosen vertex; moreover
the spectral decomposition \eqref{specdecomp} is determined by the
Kesten spectral measure at a single vertex. In particular, this
happens for Cayley graphs. In the case of the Cayley graph of
$\Ga$, one can alternatively use von Neumann traces.  The von Neumann trace on
the von Neumann algebra generated by the left regular
representation of $\Ga$ is given by
\begin{equation}\label{defineVNtrace}
\mathrm {tr}(T) = \langle T\delta_1,\delta_1\rangle;
\end{equation}
that is, it is the ``coefficient" in $T$ of the identity element
(this is literally true for elements of finite support, i.e.\
elements of $\R \Ga$).  Then $\mu^\Ga = \mathrm{tr}(E(B))$
coincides with the Kesten spectral measure (at, for example, the
identity), as one easily sees by comparing
\eqref{definespectmeasures} and \eqref{defineVNtrace}.

\subsection{Ihara zeta functions and Kesten-von Neumann-Serre spectral measures}
The following is from \cite{GZgraphs,GZihara}.  If
$X=(V,E)$ is a graph, one can define the path metric on $V$ by
setting $d(v_1,v_2)$ to be the number of edges in the shortest
path connecting $v_1$ and $v_2$.  Denote by $B_X(v,r)$ the open
ball of radius $r$ in $V$ around $v$.

Fix a positive integer $k$ as before.   The space of (isomorphism
classes of) pointed $k$-regular graphs $(X,v)$ becomes an
ultrametric compact totally disconnected space
\cite{GZgraphs,GZihara} by taking
$$d((X_1,v_1),(X_2,v_2))=\frac{1}{n+1}$$ where $n$
is the largest integer for which  $B_{X_1}(v_1,n)$ is pointedly
isometric to $B_{X_2}(v_2,n)$.

The case of interest to us is the following: $\Ga$ is a discrete
group with finite symmetric generating set $S$ and $P\leq \Ga$.
Suppose $P_n\leq \Ga$, $n\in \N$ are finite index subgroups with
$\bigcap P_n = P$ (so $P$ is a closed subgroup in the profinite
topology on $\Ga$). Then one easily sees that the graphs $X_n=\Sch
\Ga {P_n} S$ converge to the graph $X=\Sch \Ga P S$ (see
\cite{GZgraphs,GZihara}) where the coset of the identity is taken
as the base point.

A sequence of probability measures $\mu_n$ on a
measure space is said to converge \emph{weakly} to a measure $\mu$
if, for all measurable subsets $B$, $\mu_n(B)\to \mu(B)$
\cite{Billingsley}. Suppose $(X_n,v_n)$ is a sequence of pointed
finite graphs converging to $(X,v)$, and let $N>0$ be given. Then
$B_X(v,N)$ is pointedly isometric to $B_{X_n}(v_n,N)$ for $n$
sufficiently large.  From this it can be deduced using the method
of moments that $\mu^{X_n}_{v_n}\to \mu^X_{v}$ weakly
\cite{GZgraphs,GZihara}.

To motivate the so-called Kesten-von Neumann-Serre spectral
measures \cite{GZihara} we recall the definition of the
\emph{Ihara zeta function} $\zeta_X$ of a finite $k$-regular graph
$X$ \cite{Ihara}.  It is the power series
$$\zeta_X(t) = \prod _{[C]} (1-t^{|C|})\inv$$
where $[C]$ runs over the equivalence classes of primitive,
cyclically reduced closed paths in $X$ and $|C|$ denotes the
length of $C$.  It is known \cite{Ihara,GZihara} that $$\ln
\zeta_X(t) = \sum_{r=1}^{\infty} \frac{c_r}{r}t^r$$ where $c_r$ is
the number of cyclically reduced closed paths of length $r$ in
$X$. We remark that $\zeta_X$ can be viewed as a discrete analogue
of the Riemann zeta function, and satisfies the Riemann hypothesis
if and only if $X$ is a Ramanujan graph \cite{GZihara,Ramanujan}.

Let $M$ be the Markov operator associated to $X$.
 Recalling that $kM$ is the incidence matrix of the graph $X$, the
 results of \cite{Ihara} show that $$\zeta_X(t) =
 (1-t^2)^{-\frac{1}{2}(k-2)|E(X)|}\mathrm{det}(1-tkM +
 (k-1)t^2)\inv$$ where $E(X)$ denotes the edge set of $X$.

Grigorchuk and \.{Z}uk extended this to infinite graphs that are
limits of sequences of finite graphs.  Let $\{X_n\}$ be a sequence
of finite $k$-regular graphs with associated Markov operators
$M_n$.  Denote by $V_n$ the vertex set of $X_n$. We use
$\mathrm{tr}(A)$ to denote the trace of a matrix $A$. Define
\begin{equation}\label{averagemeasure}
\mu_n = \frac{1}{|V_n|} \sum_{v\in V_n} \mu^{X_n}_{v}=
\sum_{\lambda\in
  \Sp{M_n}}\frac{\mathrm{tr}(E_{\lambda})}{|V_n|}\delta_{\lambda} =
\sum_{\lambda\in \Sp {M_n}} \frac{\#_n(\lambda)}{|V_n|}
\delta_{\lambda}
\end{equation}
where $\#_n(\lambda)$ denotes the multiplicity of $\lambda$ as an
eigenvalue of $M_n$, where $E_{\lambda}$ is the projection to the
$\lambda$-eigenspace of $M_n$ and where $\delta_{\lambda}$ is the
Dirac measure at $\lambda$. The second equality of
\eqref{averagemeasure} follows from \eqref{finitedecomp}. The
probability measure $\mu_n$ counts the frequency of the
eigenvalues weighted by their multiplicities.

 Following Serre \cite{Serre} (see also \cite{GZihara}) the
eigenvalues of the Markov operators $M_n$ are said to be
\emph{equidistributed} with respect to a measure $\mu$ having
support on $[-1,1]$ if $\mu_n\to \mu$ weakly.

Suppose now that we are in the situation above with $X=\Sch \Ga P
S$ and $X_n = \Sch \Ga {P_n} S$.  It is shown in \cite{GZihara}
that the eigenvalues of the $M_n$ are equidistributed with respect
to some measure $\mu$, which Grigorchuk and \.{Z}uk call the
\emph{Kesten-von Neumann-Serre} (KNS) spectral measure of $X$ with
respect to the approximating sequence $X_n$.

Serre \cite{Serre} proved that the eigenvalues of $M_n$ are
equidistributed with respect to some measure if and only if the
sequence $\zeta_{X_n}^{\frac{1}{|E(X_n)|}}$ converges in
$\R[\!{[}t]\!{]}$ with the topology of pointwise convergence of
the coefficients. Grigorchuk and \.{Z}uk then defined the zeta
function of $X$ (in the above context) with respect to the
approximating sequence $X_n$ to be given by $$\ln \zeta_X(t) =
\lim _{n\to \infty} \frac{1}{|E(X_n)|} \ln \zeta_{X_n}(T).$$  They
showed \cite{GZihara} that $\ln \zeta_X(t)$ has radius of
convergence at least $\frac{1}{k-1}$ and that, for
\mbox{$|t|<\frac{1}{k-1}$},
\begin{equation}\label{zetaintegralform}
\ln \zeta_X(t) = -\frac{k-2}{2}\ln (1-t^2) - \int_{-1}^{1} \ln
(1-tk\lambda + (k-1)t^2) \mathrm{d} \mu\end{equation} where $\mu$
is the KNS spectral measure.  This serves as a motivation for
computing the measure $\mu$. Conversely, it is shown in
\cite{GZihara} that $\zeta_X$ determines the moments of $\mu$ and
hence determines $\mu$ itself.

They also define in \cite{GZihara} the zeta function of $\Ga$
(with respect to the generators $S$) by
$$\ln \zeta_\Ga(t) =-\frac{k-2}{2}\ln (1-t^2)
- \mathrm {tr} \ln (1-tkM + (k-1)t^2)$$ where $\mathrm{tr}$ is, as
above, the von Neumann trace.  Notice that this is the same
formula as \eqref{zetaintegralform}, but with $\mu$ replaced by
the Kesten spectral measure for the random walk on the Cayley
graph of $\Ga$.


\section{Automata groups and Cayley machines}\label{sec:preliminaries}

In this section, we introduce finite automata and the groups they generate.

A finite (Mealy) automaton \cite{Arbib,Eilenberg}
$\Ac$ is a $4$-tuple $(Q,A,\delta,\lambda)$ where $Q$ is a finite
set of states, $A$ is a finite alphabet, $\delta:Q\times A\to Q$
is the transition function and $\lambda: Q\times A\to A$ is the
output function. One writes $qa$ for $\delta(q,a)$ and $q\circ a$
for $\lambda (q,a)$. These functions extend to the free monoid
$A^*$ by
\begin{gather} q(au) = (qa)u \\\label{action} q\circ (au) =
(q\circ a)(qa)\circ u.
\end{gather} We use $\Ac_q$ to denote the
\textit{initial automaton} $\Ac$ with designated start state $q$.
There is a function $\Ac_q:A^*\to A^*$ given by $w\mapsto q\circ
w$. This function is length preserving and extends continuously
\cite{GNS} to the set of right infinite words $A^{\omega}$ via the
formula
\begin{equation}\label{actiononcantor}
\Ac_q(a_0a_1\cdots) = \lim_{n \to \infty} \Ac_q(a_0\cdots a_n)
\end{equation}
 where $A^{\omega}$
is given the product topology, making it homeomorphic to a Cantor
set. If, for each $q$, the state function $\lambda_q:A\to A$ given
by $\lambda_q(a) = q\circ a$ is a permutation, then $\Ac_q$ is an
isometry of $A^{\omega}$ for the metric $d(u,v) = 1/(n+1)$ where
$n$ is the length of the longest common prefix of $u$ and $v$
\cite{GNS}. In this case the automaton is called
\emph{invertible}.   We shall assume here (unless otherwise
stated) that all automata are invertible. Let $\Ga=\GG (\mathcal
A)$ be the group generated by the $\Ac_q$ with $q\in Q$.

If we let $T$ be the Cayley tree of $A^*$, then $\Ga$ acts on the
left of $T$ by rooted tree automorphisms of $T$ \cite{GNS,branch}
via the action \eqref{action}.  The induced action on the boundary
$\partial T$ (the space of infinite directed paths from the root)
is just the action \eqref{actiononcantor} of $\Ga$ on
$A^{\omega}$.

The automorphism group $\mathrm{Aut}(T)$ is the iterated
(permutational) wreath product of countably many copies of the
left permutation group $(S_{A},A)$ \cite{GNS,branch,Rhodestrees},
where $S_A$ denotes the symmetric group on $A$.  In this paper,
our notation will be such that the wreath product of left
permutation groups has a natural projection to its leftmost
factor; this is in contrast to the case of restricted wreath
products of abstract groups where our notation is such that there
is a projection to the rightmost factor. For a group $\Ga=\GG
(\Ac)$ generated by an automaton over $A$, one has an embedding
\begin{equation}\label{wreathembedding1}
(\Ga,A^{\omega}) \hookrightarrow (S_{|A|},A)\wr (\Ga,A^{\omega}).
\end{equation}
 The
maps sends $\Ac_q$ to the element with wreath product coordinates:
\begin{equation}\label{wreathcoords}
\Ac_q = \lambda_q(\Ac_{qa_1},\ldots,\Ac_{qa_n})
\end{equation}
where $A=\{a_1,\ldots,a_n\}$.  See \cite{GNS,branch,reset} for
more details.

The action of $\Ga$ on $T$ is called \emph{spherically transitive}
if $\Ga$ acts transitively on each level of the tree.  Here the
$k^{th}$ level of $T$ is the set of all vertices corresponding to
words of length $k$. It will be convenient to denote by $T_k$ the
finite rooted tree obtained by pruning the levels after level $k$.
Denote by $\mathrm{St}_{\Ga}(k)$ the set of all elements of $\Ga$
that fix each vertex of level $k$; it is a finite index normal
subgroup, being the kernel of the projection $\Ga\to
\mathrm{Aut}(T_k)$. Notice that $\{1\}=\bigcap_{k=0}^{\infty}
\mathrm{St}_{\Ga}(k)$ (and so $\Ga$ is residually finite).

In this paper, we shall be particularly concerned with one class of
examples. Let $G$
be a non-trivial finite group. By the \textit{Cayley machine}
$\mathcal C(G)$ of $G$ we mean the automaton with state set and
alphabet $G$. Both the transition and the output functions are the
multiplication of the group.  So in state $g_0$ on input $g$ the
machine goes to state $g_0g$ and outputs $g_0g$.  You can view
$\mathcal C(G)$ as the Cayley graph of $G$ with respect to the
generators $G$ where output is the next state.  The state function
$\lambda_g$ is just left translation by $g$ and hence a
permutation, so $\CC (G)$ is invertible. Cayley machines for
semigroups are an important part of Krohn-Rhodes theory
\cite{Arbib,KR}.  The study of the automata group of the Cayley
machine of a finite group was initiated by the second and third
authors \cite{reset}.

An automaton is called a \emph{reset automaton} if, for each $a\in
A$, $|Qa|=1$; that is, each input resets the automaton to a single
state.  The second and third authors showed that the inverse of a
state $\CC(G)_g$ is computed by the corresponding state of the
reset automaton $\mathcal A(G)$ with states $G$ and input alphabet
$G$, where in state $g_0$ on input $g$, the automaton goes to
state $g$ and outputs $g_0\inv g$.  Therefore $\GG(\CC(G))=
\GG(\Ac(G))$. In wreath product coordinates
\begin{equation}\label{wreathcoords2}
\Ac(G)_g = g\inv(\Ac_{g_1},\ldots,\Ac_{g_n})
\end{equation}
where $G=\{g_1,\ldots,g_n\}$. Hence, in our situation, there is an
embedding
\begin{equation}\label{wreathembedding}
\GG(\Ac(G))\hookrightarrow (G,G)\wr (\GG (\Ac (G)),G^{\omega}).
\end{equation}
Moreover, the action of $\GG(\Ac(G))$ on the Cayley tree of $G^*$
is spherically transitive \cite{reset}.

 Let $x= \Ac(G)_1$.  Notice that
$$x\Ac(G)\inv_g=x\CC(G)_g =g(1,\ldots,1),$$ so we can identify $G$
with a subgroup of $\GG(\Ac(G))$ via $g\leftrightarrow
x\Ac(G)\inv_g$. Let
\begin{equation}\label{defineN}
N=\langle x^nGx^{-n}\mid n\in \Z\rangle.
\end{equation}
It is shown in \cite{reset} that $x$ has infinite order, $N$ is a
locally finite group and $\GG(\Ac(G)) = N\rtimes\langle x\rangle$.
If $G$ is abelian, it is shown \cite{reset} that
$$\GG(\Ac(G)) = G\rwr \Z = \left(\bigoplus_{\Z} G\right)\rtimes \Z$$ where in
the latter semidirect product, $\Z$ acts by the shift.   In
particular, $\GG(\Ac(\Z/n\Z))$ is the lamplighter-type group
$\Z/n\Z \rwr \Z$.

The \emph{depth} of an element
$\ga\in \Ga$ is the least integer $n$ (if such exists) so that
$\ga$ only changes the first $n$ letters of a word.  An automorphism
of finite depth is often called finitary.
An important role is also played by the subgroup
\begin{equation}\label{defineN0}
N_0 = \langle x^nGx^{-n}\mid n \geq 0\rangle.
\end{equation}
It is shown in \cite{reset} that $x^ngx^{-n}$ has depth $n+1$ and
so $N_0$ consists of finitary automorphisms.   We recall also from
\cite{reset} that the elements of the form $\Ac(G)_g \in \Ga$ with
$g \in G$ generate a free subsemigroup of $\Ga$ (this holds for a
large class of semigroups generated by invertible reset automata
\cite{reset}).

\section{Dynamics, free actions and spectral measures}

Let $A$ be a finite alphabet and $T$ be the Cayley tree of $A^*$.
Let $\nu$ be the product measure on $\partial T=A^{\omega}$; so
for each $u \in A^*$, the cylinder set $u A^{\omega}$ is given the
measure $1/|A|^{|u|}$.  Let $\p:T\to T$ be a morphism of
trees preserving the distance from the root (and hence the levels);
equivalently 
$\p:A^*\to A^*$ is a function computed by a possibly infinite
automaton~\cite{GNS}; such a map is called in \cite{Rhodestrees} an
\emph{elliptic contraction}. We also use $\p:\partial T\to \partial T$ to
denote the induced morphism on the boundary (which is a
contraction~\cite{GNS}). Define $$\mathrm{Fix}(\p)=\{w\in
\partial T\mid \p(w)=w\}.$$  It is a closed subset of $\partial
T$.  Also define, for $k\geq 0$, $$\mathrm{Fix}_k(\p) = \{w\in A^k\mid
\p(w)=w\}.$$

\begin{proposition}\label{measurefixedpoints}
Let $T$ be the Cayley tree of $A^*$ and $\p:T\to T$ be an elliptic
contraction.  Then  
\begin{equation}\label{approximatefix}
\nu(\mathrm{Fix}(\p)) = \lim_{k\to\infty}
\frac{1}{|A|^k}{|\mathrm{Fix}_k(\p)|}.
\end{equation}
In particular, \eqref{approximatefix} holds for rooted automorphisms of
$T$.\qed
\end{proposition}
\begin{proof}
Set $F_k =\bigcup \{wA^{\omega}\mid w\in \mathrm{Fix}_k(\p)\}$.
Then $\{F_k\}_{k\geq 0}$ is a decreasing sequence of measurable sets and
$$\mathrm{Fix}(\ga) = \bigcap_{k=0}^{\infty} F_k.$$
Hence, since $\nu$ is a finite measure,
$$\nu(\mathrm{Fix}(\p)) = \lim_{k\to\infty} \nu (F_k) =\lim_{k\to\infty}
\frac{1}{|A|^k}|\mathrm{Fix}_k(\p)|,$$ thereby establishing the
proposition.
\end{proof}

Recall that a subset of a topological space is said to be
\emph{nowhere dense} if it does not contain a non-empty open
subset. In general, measure zero is very far from nowhere dense --
for instance the irrational numbers in the interval $[0,1]$ are
nowhere dense but have full measure.   However, for the case of a
contraction of the Cantor set computed by a finite state automata we
can prove the 
following (somewhat surprising) result, which is inspired by a
much more complicated argument from \cite{GZlamp} for the case of
the lamplighter group.

\begin{theorem}\label{nowheredensethm}
Let $\p:A^{\omega}\to A^{\omega}$ be a function computed by a
finite state automaton.  Then $\mathrm{Fix}(\p)$ has measure zero
if and only if it is nowhere dense.
\end{theorem}
\begin{proof}
Since each non-empty open subset of $A^{\omega}$ has positive
measure, it suffices to show that if $\mathrm{Fix}(\p)$ is nowhere
dense, then it has measure zero.  Let $\Ac= (Q,A,\delta,\lambda)$
be a finite state automaton such that $\p=\Ac_{q_0}$, $q_0\in Q$.
Let $Q'$ be the set of $q\in Q$ such that $\Ac_q$ is not the
identity. The hypothesis on $\p$ then implies that if $u\in A^*$
is fixed by $\p$, then $q_0\cdot u\in Q'$ -- if this were not the
case $\p$ would fix $uA^{\omega}$. Since $Q$ is finite, we can
find an integer $p> 0$ such that, for each $q\in Q'$, $\Ac_q$ does
not fix some element of $A^p$. Let $n=|A|$.  We claim:
\begin{equation}\label{numfix}
|\mathrm{Fix}_{pk}(\p)| \leq
(n^p-1)^k
\end{equation}
 for all $k\geq 0$.  We proceed by induction.  The result is clear for
 $k=0$.  Suppose that \eqref{numfix}
holds for $k-1$, with $k\geq 1$, and consider the $n^{pk}$ words
of length $pk$. By the inductive hypothesis, at most $(n^p -
1)^{k-1}$ of the possible prefixes of length $p(k-1)$ of such
words are fixed by $\p$. Now let $u\in A^{p(k-1)}$ be such a
prefix that is fixed by $\p$. Then $q_0\cdot u\in Q'$ by the
hypothesis
 on $\p$.  By choice of $p$, $\Ac_{q_0\cdot u}$ does not fix all words of length
 $p$,  so there are at most $n^p -1$ words of length $p$ that are
 fixed by $\Ac_{q_0\cdot u}$. Hence, there
are at most $(n^p - 1)^{k-1}\cdot (n^p - 1) = (n^p-1)^k$ words
 from $A^{pk}$ that are fixed by $\p$, establishing \eqref{numfix}.
Thus $$\frac{1}{n^{pk}}|\mathrm{Fix}_{pk}(\p)| \leq
\left(1-\frac{1}{n^p}\right)^k.$$  Since the right hand side of
the above equation tends to $0$ as $k\to \infty$, the theorem
follows by an application of Proposition \ref{measurefixedpoints}.
\end{proof}

Let $\Ga$ be a group acting spherically transitively by rooted
automorphisms of $T$.   Then $\Ga$ acts ergodically on the measure
space $(\partial T,\nu)$ by measure preserving transformations and
topologically transitively on $\partial T$ by isometries
\cite{GZlamp,GNS}. Let $\pi:\Ga\to \mathcal B(L^2(\partial
T,\nu))$ be the associated representation, $S$ be a finite
symmetric generating set for $\Gamma$ and $H_{\pi}$ be the
associated Hecke-type operator. Let $\pi_k:\Ga\to \mathcal
B(\ell^2(A^k))$ be the permutation representation of $\Ga$ on the
$k^{th}$ level of the tree. Identifying the elements of the $A^k$
with the characteristic functions of the associated cylinder sets,
$\pi_k$ can be viewed as a subrepresentation of $\pi$ and
$\pi_{k-1}$ as subrepresentation of $\pi_k$.  If, for $k\geq 1$,
$\pi_k'$ is the orthogonal complement of $\pi_{k-1}$ in $\pi_k$
and $\pi_0'=\pi_0$, then \cite{Hecke} $\pi =
\bigoplus_{k=0}^{\infty} \pi_k'$ from which one
obtains~\cite{Hecke}
\begin{equation}\label{pispec}
\Sp {H_{\pi}} = \ol{\bigcup_{n\geq 0} \Sp {H_{\pi_k}}}
\end{equation}

Recall \cite{Hall} that if $\alpha$ is a finite dimensional
representation of $\Ga$, then the associated character
$\chi_{\alpha}$ is given by $\chi_{\alpha}(\ga) =
\mathrm{tr}(\alpha(\ga))$.  If $\alpha$ is a permutation
representation, then $\chi_{\alpha}$ is the \emph{fixed-point}
character; it simply counts the number of fixed-points of each
element of $\Ga$.   It then seems natural to define the
\emph{fixed-point character} of $\pi$ by
$$\chi_{\pi}(\ga) = \nu(\mathrm{Fix}(\ga)).$$
An immediate consequence of Proposition \ref{measurefixedpoints} is
\begin{equation}\label{existenceoflimitcharacter}
 \chi_{\pi}(\ga) = \lim_{k\to\infty}
\frac{1}{n^k}\chi_{\pi_k}(\ga).
\end{equation}

An action of a group $\Ga$ on a measure space is said to be
\emph{free in the sense of ergodic theory} if, for all $1\neq
\ga\in \Ga$, $\mathrm{Fix}(\ga)$ has measure zero. We say that the
action is \emph{free in the sense of Baire category} if, for all
$1\neq \ga\in \Ga$, $\mathrm{Fix}(\ga)$ is nowhere dense.
  Notice that in either
of these two cases there is an infinite path $w\in \partial T$ with
trivial stabilizer.  Indeed, $\bigcup_{1\neq \ga\in \Ga}
\mathrm{Fix}(\ga)$ cannot be all of $\partial T$: in the first
case this set has measure zero; in the latter it is a countable
union of closed nowhere dense sets and hence cannot be all of
$\partial T$ by the Baire category theorem.
Theorem \ref{nowheredensethm} then has the following interpretation.

\begin{theorem}\label{tracecriterion}
Let $\Ga$ be a group acting on the Cayley tree $T$ of $A^*$ (with $A$
finite) by rooted automorphisms computed by finite state automata.
Then the following are equivalent:
\begin{enumerate}
\item The action of $\Ga$ on $(\partial T,\nu)$ is free in the sense of ergodic theory;
\item The action of $\Ga$ on $\partial T$ is free in the sense of Baire
category;
\end{enumerate}\qed
\end{theorem}

In light of the above result, it seems natural to say that the
action of a group $\Ga$ generated by a finite state automaton is
\emph{free} if the equivalent conditions of Theorem
\ref{tracecriterion} hold. The only non-trivial examples of free
actions of automata groups in the literature, so far as we know,
are the Cayley machines of finite abelian groups \cite{reset},
including the automaton of the lamplighter group considered in
\cite{GZlamp}. For Cayley machines of non-abelian groups,
fixed-point sets can have non-empty interior \cite{reset}.  It
would be interesting to find more examples.

Returning to the spectra of Hecke operators, let us fix $w\in
\partial T$ and set $w_k$ to be the prefix of $w$ of length $k$.
Let $P = \mathrm{St}_{\Ga}(w)$ and $P_k = \mathrm{St}_{\Ga}(w_k)$.
The subgroup $P$ is called a \emph{parabolic subgroup} \cite{GNS}.
Then, for all $k\geq 0$, $P_k$ is of finite index in $\Ga$ and
$P=\bigcap_k P_k$. Moreover, since $\Ga$ acts spherically
transitively, $\pi_k$ is equivalent to the quasi-regular
representation $\lambda_{\Ga/P_k}$ and $H_{\pi_k} = M_k$ the
Markov operator on \Sch {\Ga} {P_k} S; if $M$ is the Markov
operator on \Sch {\Ga} P S, then $H_{\lambda_{\Ga/P}} = M$. If, in
addition, $\Ga$ is amenable, then \cite{Hecke}[Theorem 3.6] shows
that
\begin{equation}\label{otherspec}
\Sp {H_{\pi}} = \ol {\bigcup_{n\geq 0} \Sp {H_{\pi_k}}}=\Sp
{H_{\lambda_{\Ga/P}}} \subseteq \Sp {H_{\lambda_\Ga}}\subseteq
[-1,1]
\end{equation}

We now relate the fixed-point character to the moments of the KNS
spectral measure $\mu$ associated to $\Sch {\Ga} P S$ with respect
to the approximating sequence $\Sch {\Ga}  {P_k} S$.  If $w\in
S^*$, we use $[w]$ to denote the image of $w$ in $G$.  Notice
that, for $m\geq 0$,
\begin{equation}\label{powerofM}
M^m = \frac{1}{|S|^m}\sum_{w\in S^m} [w].
\end{equation}
 The moments of $\mu_k$ (see definition
\eqref{averagemeasure}) are then given by:
\begin{equation*}
\mu_k^{(m)} = E_{\mu_k}[\lambda^m] = \sum_{\lambda\in \Sp
{M_k}}\lambda^m \frac{\#(\lambda)}{n^k} =
\frac{1}{n^k}\mathrm{tr}(M_k^m) = \frac{1}{|S|^m}\sum_{w\in S^m}
\frac{1}{n^k}\chi_{\pi_k}([w])
\end{equation*}
where the last equality holds from \eqref{powerofM} and the fact that
$M_k=\pi_k(M)$.  
Notice $0\leq \mu_k^{(m)}\leq 1$, indeed it is the average over
all $w\in S^m$ of the probability of $\pi_k([w])$ fixing a vertex
on the $k^{th}$ level.  Thus the moment generating function of
$\mu_k$ is analytic on $\mathbb{R}$ for all $k\geq 0$.   Since the
$\mu_k$ converge weakly to $\mu$, it follows that, for each $m\geq
0$, $\mu_k^{(m)}\to \mu^{(m)}$ \cite{Billingsley}. We are then led by
\eqref{existenceoflimitcharacter} to the following
formula for the moments of $\mu$:
\begin{equation}\label{momentsofmu}
\mu^{(m)} = \frac{1}{|S|^m}\sum_{w\in S^m}\chi_{\pi}([w]).
\end{equation}
Thus the $m^{th}$ moment of $\mu$ is the average over all $w\in
S^m$ of the probability of $\pi([w])$ fixing an infinite path,
whence $0\leq \mu^{(m)}\leq 1$.  This implies that the moment
generating function of $\mu$ is analytic on $\mathbb{R}$.

Let us contrast this with the situation for the Kesten spectral
measure $\mu^{\Ga}$  for the random walk on the Cayley graph of
$\Gamma$.  As mentioned earlier, the moments correspond to return
probabilities.  More precisely, if the Markov operator $M$ on
$\ell^2(\Ga)$ has spectral decomposition
$$M = \int_{-1}^1 \lambda \mathrm dE(\lambda)$$ with (projection-valued)
spectral measure $E$, then
$$\mu^{\Ga}(B) = \langle \int_B\mathrm
dE(\lambda)\delta_1,\delta_1\rangle$$ for $B$ a Borel subset of
$[-1,1]$.  So the $m^{th}$ moment is given by
\begin{align*}
(\mu^{\Ga})^{(m)} &= \int_{-1}^1 \lambda ^m\mathrm d\mu^{\Ga} =
\langle \int_{-1}^1 \lambda^m\mathrm
dE(\lambda)\delta_1,\delta_1\rangle \\ &= \langle
M^m\delta_1,\delta_1\rangle = \mathrm{tr}(M^m)
\end{align*} where the last trace is the von Neumann trace, c.f.\
\eqref{defineVNtrace}. Since $M^m\in \R\Ga$, this is just the
coefficient of $1$. Notice that the return probability $p_m(1)$ is
just the fraction of words in $S^m$ representing the identity $1$
of $\Ga$. It follows from \eqref{powerofM} that $(\mu^{\Ga})^{(m)}
= p_m(1)$. From this, we may easily deduce that the moment
generating function of $\mu^{\Ga}$ is analytic on all of
$\mathbb{R}$.

Now we are in a position to compare $\mu^{(m)}$ with
$(\mu^{\Ga})^{(m)}$. Recalling that \mbox{$0\leq
\chi_{\pi}(\ga)\leq 1$} and $\chi_{\pi}(1)=1$, it is clear that
the right hand side of \eqref{momentsofmu} is $p_m(1)$ precisely
if, for each $w\in S^m$ such that $[w]\neq 1$,
$\chi_{\pi}([w])=0$. In other words, the average probability that
a word in $w\in S^*$ of length $m$ fixes an infinite path will be
the same as the probability that $w$ represents the identity if
and only if for all $w\in A^*$, either $w$ represents the
identity, or $w$ almost surely does not fix any infinite path.
Recall \cite[Theorem 30.1]{Billingsley} that if two probability
measures on $[-1,1]$ have the same moment generating function, and
this function is analytic on a neighbourhood of $0$, then the
measures are the same.   Recalling that $\chi_{\pi}(\ga) =
\nu(\mathrm{Fix}(\ga))$, we may summarize the previous discussion
in the following theorem.

\begin{theorem}\label{freevs1}
Let $\Ga$ be a group acting spherically transitively on the Cayley
tree of $A^*$ with finite symmetric generating set $S$. Let $w \in
A^{\omega}$ be an infinite path and, for all $k\geq 0$, let $w_k$
be the prefix of $w$ of length $k$.  Set $P=\mathrm{St}_{\Ga}(w)$
and $P_k=\mathrm{St}_{\Ga}(w_k)$ and let $X$ and $X_k$ be the
respective Schreier graphs with respect to a finite generating set
$S$. Then the following are equivalent:
\begin{enumerate}
\item The KNS spectral measure for $X$, with respect to the
approximating sequence $\{X_k\}$, coincides with the Kesten
spectral measure for the simple random walk on the Cayley graph of
$\Ga$
\item The action of $\Ga$ is free in the sense of ergodic theory.
\end{enumerate}
\qed
\end{theorem}

An immediate consequence of Theorems \ref{tracecriterion} and
\ref{freevs1} is:
\begin{corollary}\label{zetafunction}
Let $\Ga$ be an automata group acting freely on the boundary of
$T$.  Choose $w\in \partial T$ with trivial stabilizer and define
$P_k$, $k\geq 0$, as per Theorem \ref{freevs1}. Then
$\zeta_{\Ga}$ coincides with the limit zeta function from the
approximating sequence $\Sch \Ga {P_k} S$.\qed
\end{corollary}

\section{Calculation of the spectral measures for Cayley machines}

In this section, we present a computation of the spectra and spectral
measures associated to random walks on automata groups generated by
the Cayley machines of finite groups. This not only provides an
interesting new class of examples of spectra and spectral measures of
random walks, but also serves to illustrate the material introduced
in the previous sections.

Fix for this section a non-trivial finite group $G=\{g_1,\ldots,
g_n\}$. As standing notation we set $\Ga = \GG(\Ac(G))$.   It is
locally finite-by-infinite cyclic \cite{reset} and hence amenable
so \eqref{otherspec} applies to computing the spectrum.  In the
case that $G$ is abelian, the results of \cite{reset} show the
action of $\Ga$ is free.

Let us establish some notation. We use $\ol{g_i}$ as a shorthand
notation for the element $\Ac_{g_i} \in \Ga$. Let
\begin{equation}\label{generatorsS}
S=\{\ol{g_1},\ldots \ol{g_n},\ol{g_1}\inv,\ldots
\ol{g_n}\inv\};\end{equation} so $|S|=2n$. We fix a parabolic
subgroup $P$. If $G$ is abelian, we choose $P$ to be trivial. If
$G$ is non-abelian, we can choose $P$ to be a locally finite group
\cite{reset}.   Let $w$ be an infinite path with stabilizer $P$
and for $k \geq 0$ let $w_k$ denote the prefix of $w$ of length
$k$. Let $P_k \leq \Gamma$ be the stabilizer of $w_k$.  Now $P_k$
has index $n^k$ in $\Gamma$ (by spherical transitivity) and
$P=\bigcap P_k$. We shall show that all the spectra considered in
equation \eqref{otherspec} are the entire interval $[-1,1]$.  We
shall also calculate the KNS spectral measure associated to the
Schreier graph $X=\Sch {\Ga} P S$ with respect to the
approximating sequence of graphs $X_k = \Sch {\Ga} {P_k} S$, in
particular obtaining the Kesten spectral measure for $\Ga$ in the
case $G$ is abelian.

\subsection{Operator recursion, wreath products and the monomial
representation}\label{oprecursion}

We begin by recalling some standard facts about the matrix
representation of wreath products of permutations groups
\cite{Hall,Arbib}. Let $(H,X)$ and $(K,Y)$ be left permutation
groups and let $(W,X\times Y)= (H,X)\wr (K,Y)$. The associated
\emph{monomial representation} \cite{Hall,Arbib} is described as
follows.  Let $(h,f)\in W=H\times K^X$. The \emph{monomial matrix}
for $(h,f)$ is obtained from the $|X|\times |X|$ permutation
matrix for $h$ by replacing the $1$ in column $i\in X$ by $f(i)$.
The action on $X\times Y$ is recovered by considering column
vectors of size $|X|$ with entries in $Y$.  To linearize this
representation, one has to replace column vectors of size $|X|$
with entries in $Y$ by size $|X|$ block column vectors with
entries size $|Y|$ column vectors. Then the matrix representation
associated to $(W,X\times Y)$ takes the $|X|\times |X|$ monomial
matrix corresponding to an element $(h,f)$ and replaces each entry
from $k\in K$ with the $|Y|\times |Y|$ permutation matrix
associated to $k$ from the permutation representation $(K,Y)$ and
replaces the zeroes by the $|Y|\times |Y|$ zero matrix. Thus the
matrices for $(W,X\times Y)$ are block monomial.

Observe that the wreath product coordinates \eqref{wreathcoords2},
restricted to $T_{k+1}$, show
$$(\pi_{k+1}(\Ga),G^{k+1})\leq (G,G)\wr (\pi_k(\Ga),G^k).$$ On the
automata generators, the embedding is given by
\begin{equation}\label{oprec1}
\olg{i}{k+1}\longmapsto g_i\inv(\olg{1}{k},\ldots, \olg{n}{k})
\end{equation}
This lets us construct inductively the matrices for the
representations $\pi_k$ using the monomial representation; this
procedure is called \emph{operator recursion} \cite{GZlamp,Hecke}.
For $i=1,\ldots, n$, set $\olg{i}{0} =1$.  Then, for $k\geq 0$,
the matrix for $\olg{i}{k+1}$ is obtained by taking the
permutation matrix from the left regular representation of $G$
corresponding to $g_i\inv$ and replacing the $1$ in column $j$ by
$\olg{j}{k}$ and the zeroes by the $n^k\times n^k$ zero matrix. So
$\olg{i}{k+1}$ is an $n\times n$ block monomial matrix with blocks
of size $n^{k}$. There is exactly one non-zero block in column
$j$, namely in row $l$, where $g_i\inv g_j = g_l$; this block is
$\olg{j}{k}$. Since all these matrices are permutation matrices,
the inverse of a matrix $\olg{i}{k+1}$ is simply the transpose
$\olgt{i}{k+1}$.

For example, if $G=\Z/2\Z=\{a,b\}$ where $a$ is the identity and
$b$ the non-trivial element, then (c.f.\ \cite{GZlamp})
\begin{equation*}
\pi_{k+1}(\ol{a}) = \begin{pmatrix} \pi_k(\ol{a}) & 0\\ 0
&\pi_k(\ol{b})\end{pmatrix},\quad \pi_{k+1}(\ol{b})
=\begin{pmatrix} 0 & \pi_k(\ol{b}) \\ \pi_k(\ol{a}) &
0\end{pmatrix}.
\end{equation*}
If $G=\Z/3\Z=\{a,b,c\}$ where $a$ is the identity, then
\begin{gather*}
\pi_{k+1}(\ol{a}) = \begin{pmatrix} \pi_k(\ol{a}) & 0 & 0\\ 0 &\pi_k(\ol{b})& 0\\
0&0&\pi_k(\ol{c})\end{pmatrix},\ \pi_{k+1}(\ol{b}) =\begin{pmatrix} 0 & \pi_k(\ol{b}) & 0\\
0& 0& \pi_k(\ol{c})\\ \pi_k(\ol{a}) &0 &0\end{pmatrix}\\
\pi_{k+1}(\ol{c})  =\begin{pmatrix} 0& 0& \pi_k(\ol{c})\\
\pi_k(\ol{a}) &0 &0\\0 & \pi_k(\ol{b}) & 0\\\end{pmatrix}.
\end{gather*}

Given square matrices $A$ and $B$ we write $A \otimes B$ for their
tensor product, given in block matrix form:
$$A\otimes B=[A_{ij} B]_{i,j = 1}^n.$$
Let $T$ be the $n \times n$ matrix defined by
$$T_{ij} = \begin{cases} 1 \text{ if } i \neq j \\
                         0 \text{ if } i = j \end{cases}$$
and set
$$S_0 = n-1, \quad S_k = T \otimes I_{n^{k}},\ k>0.$$
Notice that $T$ is the sum of the permutation matrices
corresponding to the non-identity elements of $G$ under the left
regular representation; indeed for each $i\neq j$ there is a
unique non-identity permutation of $G$ taking $i$ to $j$.

We shall need the following lemma later.
\begin{lemma}\label{addupnonidentity}
For all $k \in \N$, we have
$$\sum_{i,j = 1, i \neq j}^{n} \olg{i}{k} \olgt{j}{k} = nS_k$$
\end{lemma}
\begin{proof}
For $k=0$ this is clear, since, for each $i=1,\ldots, n$, we add
up $1$ exactly $n-1$ times.  For $k>0$, observe that
\eqref{oprec1} easily implies that in wreath product coordinates
$$\olg{i}{k} \olg{j}{k}\inv = g_i\inv g_j(1,\ldots,1).$$  So in
matrix form, we obtain $A\otimes I_{n^{k}}$ where $A$ is the
matrix for $g_i\inv g_j$ in the regular representation. If we fix
$i$ and sum over $j\neq i$ we get the sum of all permutation
matrices from the regular representation of $G$ except the
identity, that is $T$, tensored with $I_{n^{k}}$.  But this is
precisely $S_k$. Since we have $n$ choices for $i$, we obtain the
formula of the lemma.
\end{proof}

\subsection{Calculation of the characteristic polynomial}
It will be convenient to not have to always divide by $2n$ so set
$A_k=2nM_k$;  $A_k$ is the incidence matrix of the Schreier graph
$\Sch {\Ga} {P_k} S$. More explicitly,
$$A_k = \sum_{i = 1}^n (\olg{i}{k} + \olgt{i}{k}).$$

Our first objective is to calculate the spectrum of the matrix
$A_k$. To this end, define a function of two variables by
$$\Phi_k(\lambda, \mu) = \left| A_k - \lambda I_{n^k} - \mu S_k \right|$$
so that $\Phi_k(\lambda, 0)$ is the characteristic polynomial of
$A_k$. Our objective, then, is to find the roots of
$\Phi_k(\lambda, 0)$. To facilitate this, we shall obtain a
recursive formula for $\Phi_{k+1}$ in terms of $\Phi_k$. Of course
\begin{equation}\label{Phi0}
\Phi_{0}(\lambda,\mu) = 2n -\lambda -(n-1)\mu
\end{equation}
The additional term $\mu S_k$ serves as a garbage collecting term;
it arises when one tries to express the characteristic polynomial
of $A_{k+1}$ in terms of the characteristic polynomial of $A_k$.

The matrices $\olg{i}{k}$ and $\olgt{i}{k}$, $i=1,\ldots, n$, are
block monomial matrices coming from the image of $g_i\inv$,
respectively, $g_i$, under the left regular representation of $G$.
Since the sum of the permutation matrices from the left regular
representation of $G$ is the matrix of all ones (see the
discussion above concerning $T$, but now add in the identity
matrix), we see that

\begin{equation*}\label{lemma_Sk_product}
A_{k+1} = \left( \begin{matrix}
\olg{1}{k} + \olgt{1}{k} & \olg{2}{k} + \olgt{1}{k} & \dots & \olg{n}{k} + \olgt{1}{k} \\
\olg{1}{k} + \olgt{2}{k} & \olg{2}{k} + \olgt{2}{k} & \dots & \olg{n}{k} + \olgt{2}{k} \\
\vdots & \vdots &\ddots& \vdots \\
\olg{1}{k} + \olgt{n}{k} & \olg{2}{k} + \olgt{n}{k} & \dots &
\olg{n}{k} + \olgt{n}{k}
\end{matrix} \right)
\end{equation*}
and so $A_{k+1}-\lambda I_{n^{k+1}}-\mu S_{k+1}$ is given by
{\small
$$\begin{pmatrix}
\olg{1}{k} + \olgt{1}{k} - \lambda I & \olg{2}{k} +
\olgt{1}{k} - \mu I & \dots
& \olg{n}{k} + \olgt{1}{k} - \mu I \\
\olg{1}{k} + \olgt{2}{k} - \mu I & \olg{2}{k} + \olgt{2}{k}
- \lambda I & \dots
& \olg{n}{k} + \olgt{2}{k} - \mu I \\
\vdots & \vdots & \ddots& \vdots \\
\olg{1}{k} + \olgt{n}{k} - \mu I & \olg{2}{k} + \olgt{n}{k}
- \mu I & \dots & \olg{n}{k} + \olgt{n}{k} - \lambda I
\end{pmatrix}$$}
where $I$ denotes the $n^k\times n^k$ identity matrix.

We now apply some row and column operations at the block level,
designed to simplify the computation of the determinant. Applying
the operation $C_i\mapsto C_i-C_n$, for $i=1,\ldots, n-1$, yields
the matrix {\tiny
$$\begin{pmatrix}
\olg{1}{k} - \olg{n}{k} - \tsf{} & \olg{2}{k} - \olg{n}{k} & \dots
& \olg{n}{k} + \olgt{1}{k} - \mu I \\
\olg{1}{k} - \olg{n}{k} & \olg{2}{k} + \olg{n}{k} - \tsf{} & \dots
& \olg{n}{k} + \olgt{2}{k} - \mu I \\
\vdots & \vdots & \ddots& \vdots \\
\olg{1}{k} - \olg{n}{k} + (\lambda-\mu)I & \olg{2}{k} - \olg{n}{k}
+ (\lambda-\mu) I & \dots & \olg{n}{k} + \olgt{n}{k} - \lambda I
\end{pmatrix}$$}

Applying the operation $R_i\mapsto R_i-R_n$, for $i=1,\ldots,
n-1$, we obtain the matrix ($\ast$): {\Tiny
$$
\begin{pmatrix}
\tsf{-2} & \tsf{-} & \dots & \olgt{1}{k} - \olgt{n}{k}+\tsf{} \\
\tsf{-} & \tsf{-2} & \dots & \olgt{2}{k} - \olgt{n}{k}+\tsf{} \\
\vdots  & \vdots  & \ddots    & \vdots \\
\tsf{-} & \tsf{-} & \dots & \olgt{n-1}{k} - \olgt{n}{k}+\tsf{} \\
\olg{1}{k} - \olg{n}{k} +\tsf{}& \olg{2}{k} - \olg{n}{k}+\tsf{} &
\dots & \olg{n}{k}+\olgt{n}{k} - \lambda I
\end{pmatrix}$$}

To calculate the determinant of this matrix, we need some
technical results.

\begin{lemma}\label{lemondet}
Suppose that we have a block matrix
$$A=\begin{pmatrix}A_{11} & 0    &\ldots  &0 & A_{1n}\\
                 0      &A_{22}&  0      &\ldots&A_{2n}\\
                 \vdots & 0    &\ddots    & 0      &\vdots\\
                 0      &\vdots&     0    &A_{n-1,n-1} &
                 A_{n-1,n}\\
                 A_{n1} & A_{n2} &\ldots &A_{n,n-1} &
                 A_{nn}\end{pmatrix}$$
where the $A_{ij}$ are square matrices of the same size.  Moreover, suppose that
$A_{11},A_{22},\ldots,A_{n-1,n-1}$ commute with all the other
matrices. Then \begin{equation}\label{eqlemondet}|A| = \left|
A_{11}\cdots A_{nn} - \sum_{i=1}^{n-1}A_{11}\cdots
\widehat{A_{ii}}\cdots
A_{n-1,n-1}A_{ni}A_{in}\right|\end{equation} where
$\widehat{A_{ii}}$ means omit $A_{ii}$.
\end{lemma}
\begin{proof}
Since the invertible matrices are dense in the space of matrices,
we may assume without loss of generality that $A_{ii}$ is
invertible, $i=1,\ldots,n-1$.  Let $I$ be the identity matrix of
the same size as the $A_{ij}$.   Then one verifies directly that
\begin{align*}
A&=\begin{pmatrix}A_{11} & 0    &\ldots  &0 &  0\\
                 0      &A_{22}&  0      &\ldots& 0\\
                 \vdots & 0    &\ddots    & 0      &\vdots\\
                 0      &\vdots&     0    &A_{n-1,n-1} &
                 0\\
                 A_{n1} & A_{n2} &\ldots &A_{n,n-1} &
                 I\end{pmatrix}\qquad \times \\ & \qquad \qquad
                 \begin{pmatrix}I & 0    &\ldots  &0 & A_{11}\inv A_{1n}\\
                 0      &I&  0      &\vdots&A_{22}\inv A_{2n}\\
                 \vdots & 0    &\ddots    & 0      &\vdots\\
                 0      &\vdots&     0    &I &
                 A_{n-1,n-1}\inv A_{n-1,n}\\
                  0 & 0 &\ldots &0 &
                 A_{nn}-\sum_{i=1}^{n-1}A_{ni}A_{ii}\inv A_{in}\end{pmatrix}
                 \end{align*}
Using that the determinant of a block upper (lower) triangular
matrix is the product of the determinant of the diagonal blocks
and that $A_{ii}$, $i=1,\ldots,n-1$, commutes with the remaining
matrices gives \eqref{eqlemondet}.
\end{proof}

\begin{corollary}\label{cor_determinant}
Consider a block matrix
$$M = \begin{pmatrix}
2A & A & \dots & A & B_{1n} \\
A & 2A & \dots & A & B_{2n} \\
\vdots & \vdots &  & \vdots & \vdots \\
A & A & \dots & 2A & B_{n-1, n} \\
B_{n1} & B_{n2} & \dots & B_{n,n-1} & B_{nn}
\end{pmatrix}$$
where $A$ is a square matrix and $B_{1n}, \dots, B_{nn}, B_{n1},
\dots, B_{nn}$ are square matrices of the same size as $A$,
commuting with $A$.  Then $$|M|=\left|A^{n-2}\left(nAB_{nn} -
(n-1)\sum_{i-1}^{n-1}B_{ni}B_{in} + \sum_{i\neq
j}^{n-1}B_{ni}B_{jn}\right)\right|.$$
\end{corollary}
\begin{proof}

We proceed by applying the following elementary row and column
operations to the rows and columns of \textit{blocks} in $M$.
\begin{itemize}
\item[(i)] $C_i \mapsto C_i - C_{i+1} \text{ for } i = 1, \dots,
n-2$;
 \item[(ii)] $R_i \mapsto \sum_{j = 1}^i R_j\text { for } i =
n-1, \dots, 2;$ \item[(iii)] $C_{n-1} \mapsto C_{n-1} - \sum_{j =
1}^{n-2} j C_j.$
\end{itemize}
These operations leave the determinant unchanged and it is easy to
verify that they result in the matrix: {\tiny
$$\begin{pmatrix}
A & 0 & \dots & 0 & 0 & B_{1n} \\
0 & A & \dots & 0 & 0 & B_{1n} + B_{2n} \\
\vdots & \vdots & \ddots & \vdots & \vdots & \vdots \\
0 & 0 & \dots & A & 0 & \sum_{j=1}^{n-2} B_{jn} \\
0 & 0 & \dots & 0 & nA & \sum_{j=1}^{n-1} B_{jn} \\
B_{n1} - B_{n2} & B_{n2} - B_{n3} & \dots & B_{n,n-2}- B_{n,n-1} &
(n-1) B_{n,n-1} - \sum_{j=1}^{n-2} B_{nj} & B_{nn}
\end{pmatrix}$$}
Since this matrix has the same determinant as $M$, applying Lemma
\ref{lemondet} gives us that
\begin{equation}\label{messycalc1}
\begin{split}
|M|&=  nA^{n-1} B_{nn} - nA^{n-2}\sum_{i = 1}^{n-2}
\left([B_{ni} - B_{n,i+1}] \sum_{j=1}^i B_{jn}\right)\\
  & - A^{n-2} \left((n-1) B_{n,n-1} - \sum_{i=1}^{n-2} B_{ni}\right) \left( \sum_{j=1}^{n-1} B_{jn} \right).
\end{split}
\end{equation}
By telescoping we obtain
\begin{equation}\label{messycalc2}
\sum_{i = 1}^{n-2} \left([B_{ni} - B_{n,i+1}] \sum_{j=1}^i
B_{jn}\right) = \sum_{i=1}^{n-2} B_{ni}B_{in} -
B_{n,n-1}\sum_{i=1}^{n-2}B_{in}.
\end{equation}
 Substituting \eqref{messycalc2} into \eqref{messycalc1} gives
\begin{eqnarray*}
 |M| = A^{n-2}\left(nAB_{nn} - n\sum_{i=1}^{n-2}B_{ni}B_{in} +
nB_{n,n-1}\sum_{i=1}^{n-2}B_{in}  \right.\\
-\left. (n-1)B_{n,n-1}\sum_{j=1}^{n-1}B_{jn}
+\sum_{i=1}^{n-2}\sum_{j=1}^{n-1}B_{ni}B_{jn}\right)\\
=A^{n-2}\left(nAB_{nn} - n\sum_{i=1}^{n-2}B_{ni}B_{in} +
B_{n,n-1}\sum_{i=1}^{n-2} B_{in} \right.\\\left. -
(n-1)B_{n,n-1}B_{n-1,n} +
\sum_{i=1}^{n-2}\sum_{j=1}^{n-1}B_{ni}B_{jn}\right)\\
=A^{n-2}\left(nAB_{nn} - n\sum_{i=1}^{n-1}B_{ni}B_{in} +
\sum_{i=1,j=1}^{n-1}B_{ni}B_{jn}\right)\\
=A^{n-2}\left(nAB_{nn} - (n-1)\sum_{i=1}^{n-1}B_{ni}B_{in} +
\sum_{i\neq j}^{n-1}B_{ni}B_{jn}\right)
\end{eqnarray*}
as desired.
\end{proof}


Corollary \ref{cor_determinant} gives us a method of calculating
the determinant of ($\ast$).  To apply it we set
\begin{align*}
A &= (\mu-\lambda)I\\
B_{ni}&= \olg{i}{k} - \olg{n}{k} - (\mu-\lambda)I,\
i=1,\ldots,n-1\\
B_{jn} &= \olgt{j}{k}-\olgt{n}{k} - (\mu-\lambda)I,\ j=1,\ldots,
n-1\\
B_{nn} &= \olg{n}{k} + \olgt{n}{k} -\lambda I
\end{align*}
From this we obtain for $i,j=1,\ldots,n-1$
\begin{align*}
B_{ni}B_{jn} &= \olg{i}{k}\olgt{j}{k} - \olg{i}{k}\olgt{n}{k} -
\olg{n}{k}\olgt{j}{k} + I \\ & + (\mu-\lambda)^2I - (\mu
-\lambda)\left(\olg{i}{k}+\olgt{j}{k} -
\olg{n}{k}-\olgt{n}{k}\right).
\end{align*}
Notice that if $i=j$, then the first term becomes $I$.

Substituting the above values into the formula from Corollary
\ref{cor_determinant} gives
\begin{align*}
\Phi_{k+1}(\lambda,\mu) & =
(\mu-\lambda)^{(n-2)n^k}|n(\mu-\lambda)(\og{n} + \ogt{n} -\lambda I) \\
& - (n-1)\sum_{i=1}^{n-1}\left[2I - \olg{i}{k}\olgt{n}{k} -
\olg{n}{k}\olgt{i}{k}\right.  + (\mu-\lambda)^2I \\ &\left. \quad - (\mu
-\lambda)\left(\olg{i}{k}+\olgt{i}{k} -
\olg{n}{k}-\olgt{n}{k}\right)\right]\\
&+\sum_{i\neq j}^{n-1}[\olg{i}{k}\olgt{j}{k} - \olg{i}{k}\olgt{n}{k} -
\olg{n}{k}\olgt{j}{k} + I \\ & \quad + (\mu-\lambda)^2I - (\mu
-\lambda)\left(\olg{i}{k}+\olgt{j}{k} -
\olg{n}{k}-\olgt{n}{k}\right)]|\\
\end{align*}

Judicious rearranging of the various terms shows
that the above summation equals:
\begin{align*}
 \qquad &(\mu-\lambda)^{(n-2)n^k}|n(\mu-\lambda)(\og{n} + \ogt{n} -\lambda I)
\\ & -2(n-1)^2I + (n-1)\sum_{i=1}^{n-1}(\og{i}\ogt{n} + \og{n}\ogt{i})
\\ & - (n-1)^2(\mu-\lambda)^2I + (n-1)(\mu-\lambda)\sum_{i=1}^{n-1}(\og{i}+\ogt{i}) \\ &
-(n-1)^2(\mu-\lambda)(\og{n}+\ogt{n})
+ \sum_{i\neq j}^{n-1}\og{i}\ogt{j} \\ &-
(n-2)\sum_{i=1}^{n-1}(\og{i}\ogt{n} + \og{n}\ogt{i})\\ &
 + (n-1)(n-2)I +(n-1)(n-2)(\mu -\lambda)^2I \\ &- (n-2)(\mu-\lambda)\sum_{i=1}^{n-1}(\og{i}+\ogt{i})
 \\ &+ (n-1)(n-2)(\mu-\lambda)(\og{n}+\ogt{n}|
\end{align*}
Further rearrangement of the terms in attempt to obtain an expression involving $A_k$ shows that the above summation equals:
\begin{align*}
(\mu-\lambda)&^{(n-2)n^k} \times \\ & \left|(\mu-\lambda)(A_k
 -\lambda I - 
 (n-1)\mu I)-n(n-1)I  +\sum^n_{i\neq j}\og{i}\ogt{j}\right|
\end{align*}

Applying Lemma \ref{addupnonidentity} we obtain
\begin{equation}\label{recursiveformula}
\begin{split}\Phi_{k+1}&(\lambda,\mu)
 \\ &= (\mu - \lambda)^{(n-2)n^k}| (\mu - \lambda) \left( A_k -
\lambda I - (n-1) \mu I \right)\\
&\qquad\qquad\qquad\qquad\qquad
- n(n-1) I + nS_k | \\
&= (\mu - \lambda)^{(n-1)n^k}
   \left| A_k - \left(\lambda + (n-1) \mu + \frac{n(n-1)}{\mu - \lambda} \right) \right.I
   \\ &\left. \qquad\qquad\qquad\qquad\qquad + \frac{n}{\mu - \lambda} S_k \right| \\
&= (\mu - \lambda)^{(n-1)n^k}
   \Phi_k \left( \lambda + (n-1) \mu + \frac{n(n-1)}{\mu - \lambda}, \ \ - \frac{n}{\mu - \lambda}\right) \\
&= (\mu - \lambda)^{(n-1)n^k}\times \\ & \qquad
   \Phi_k \left( \frac{-\lambda^2 + (n-1) \mu^2 + (2-n) \lambda \mu + n(n-1)}{\mu - \lambda}, -\frac{n}{\mu -
   \lambda}\right)
\end{split}
\end{equation}

Next, we seek to solve the recursion and obtain an explicit
formula.

\subsection{Calculation of the eigenvalues}

Given $\lambda, \mu \in \C$, we write
$$\lambda' = \frac{-\lambda^2 + (n-1) \mu^2 + (2-n) \lambda \mu + n (n-1)}{\mu - \lambda}$$
and
$$\mu' = - \frac{n}{\mu - \lambda}$$
We define a sequence of functions $F_k(\lambda, \mu)$ inductively
by
$$F_1(\lambda, \mu) = \mu - \lambda$$
and for $k \geq 1$
$$F_{k+1}(\lambda, \mu) = F_k(\lambda', \mu').$$
\begin{lemma}\label{lemma_lambamu_facts}
Define a sequence by $(\lambda_1, \mu _1) = (\lambda, \mu)$ and
$(\lambda_{k+1}, \mu_{k+1}) = (\lambda_k', \mu_k')$ for $k \geq
1$. Then for any $k \geq 1$ we have
\begin{itemize}
\item[(i)] $\lambda_{k+1} + (n-1) \mu_{k+1} = \lambda_k + (n-1) \mu_k$;
\item[(ii)] $\mu_{k+1} - \lambda_{k+1} = - (\lambda + (n-1) \mu) - \frac{n^2}{\mu_k - \lambda_k}$;
\item[(iii)] $F_k(\lambda, \mu) = \mu_k - \lambda_k$; and
\item[(iv)] $F_{k+1}(\lambda,\mu) = -(\lambda + (n-1) \mu) - \frac{n^2}{F_k(\lambda, \mu)}.$
\end{itemize}
\end{lemma}
\begin{proof}
It is a straightforward computation to verify claims (i) and (ii)
using the definitions of
$\lambda_{k+1}$ and $\mu_{k+1}$, while claims (iii) and (iv) are
immediate consequences of (ii), together with the definitions.
\end{proof}

We remark that the rational function $f:\R^2\to \R^2$ given by
$f(\lambda,\mu) = (\lambda',\mu')$ is integrable in the sense of
\cite{GZihara}.  Namely, if $\psi:\R^2\to \R$ is given by $\psi
(\lambda,\mu) = \lambda + (n-1)\mu$ and $\alpha:\R\to \R$ is the
identity, then the previous lemma implies $\alpha\psi = \psi f$.

\begin{lemma}\label{lemma_phik_product}
For all $\mu$, $\lambda$ and all $k \geq 0$
$$\Phi_k(\lambda, \mu) = (2n - \lambda - (n-1) \mu) \prod_{i=1}^k (F_i(\lambda, \mu))^{(n-1)n^{k-i}}$$
\end{lemma}
\begin{proof}
For $k=0$, this is clear from \eqref{Phi0}. Now suppose the lemma
holds for $k\geq 0$. Then by \eqref{recursiveformula},
\begin{align*}
\Phi_{k+1}(\lambda,\mu) &=
(\mu-\lambda)^{(n-1)n^k}\Phi_k(\lambda',\mu') \\
&=(\mu-\lambda)^{(n-1)n^k}(2n-\lambda' -(n-1)\mu')\prod_{i=1}^k
(F_i(\lambda',\mu'))^{(n-1)n^{k-i}} \\
&=
(\mu-\lambda)^{(n-1)n^k}(2n-\lambda-(n-1)\mu)\prod_{i=1}^{k}(F_{i+1}(\lambda,\mu))^{(n-1)n^{k-i}}\\
&\qquad
\text{(by Lemma \ref{lemma_lambamu_facts}(i))}\\
&= (2n-\lambda -
(n-1)\mu)(F_1(\lambda,\mu))^{(n-1)n^k}\prod_{i=2}^{k+1}(F_i(\lambda,\mu))^{(n-1)n^{k+1-i}}
\end{align*}
establishing the lemma.
\end{proof}

We now want to express each $F_k$ as a rational function $P_k /
Q_k$ so that we can compute our determinant.  We define
inductively polynomials
\begin{align*}
&P_1(\lambda,\mu) = \mu-\lambda, \qquad & P_{k+1}(\lambda,\mu) & =
-(\lambda +(n-1)\mu)P_k(\lambda,\mu) - n^2 Q_k(\lambda,\mu)\\
&Q_1(\lambda,\mu) = 1, & Q_{k+1}(\lambda,\mu) &= P_k(\lambda,\mu).
\end{align*}

\begin{lemma}\label{rationalfunc}
For $k \geq 1$ we have
$$F_k(\lambda, \mu) = \frac{P_k(\lambda,\mu)}{Q_k(\lambda,\mu)}.$$
\end{lemma}
\begin{proof}
For $k = 1$ the result is clear. Now let $k \geq 1$ and
assume by induction that
$$F_k(\lambda, \mu) = \frac{P_k(\lambda,\mu)}{Q_k(\lambda,\mu)}.$$
Then Lemma \ref{lemma_lambamu_facts}(iv) tells us that
$$F_{k+1}(\lambda, \mu) \\ = - (\lambda + (n-1)\mu) - \frac{n^2}{F_k(\lambda, \mu)}$$
so we obtain \begin{align*}F_{k+1}(\lambda, \mu) &= - (\lambda +
(n-1)\mu) -
 \frac{n^2 Q_k(\lambda,\mu)}{P_k(\lambda,\mu)} \\&= \frac{-(\lambda+(n-1)\mu) P_k(\lambda,\mu)
  - n^2 Q_k(\lambda,\mu)}{P_k(\lambda,\mu)}
\\& =
 \frac{P_{k+1}(\lambda,\mu)}{Q_{k+1}(\lambda,\mu)}.\end{align*}
as required.
\end{proof}

We are primarily interested in the case $\mu=0$, so set
$P_k(\lambda) = P_k(\lambda,0)$ and $Q_k(\lambda) =
Q_k(\lambda,0)$.   Now $P_k(\lambda)$ and $Q_k(\lambda)$ satisfy:
$$\left( \begin{matrix} P_k(\lambda) \\ Q_k(\lambda) \end{matrix}
  \right) = \begin{pmatrix}  -\lambda & -n^2 \\ 1 & 0
  \end{pmatrix}\begin{pmatrix} P_{k-1}(\lambda) \\
  Q_{k-1}(\lambda)\end{pmatrix}  =
  \left( \begin{matrix} -\lambda & -n^2 \\ 1 & 0 \end{matrix} \right)^k
  \left( \begin{matrix} 1 \\ 0 \end{matrix} \right).$$

Calculating the eigenvalues and eigenvectors and then
diagonalizing, we obtain:

\begin{equation}\label{matrix_decomposition}
\begin{split}
\left( \begin{matrix} - \lambda & -n^2 \\ 1 & 0 \end{matrix} \right) = & \left( \begin{matrix} \frac{2 n^2}{- \lambda + \sqrt{\lambda^2 -
4 n^2}}
                        & \frac{2n^2}{-\lambda - \sqrt{\lambda^2 - 4 n^2}} \\ 1 & 1 \end{matrix} \right) \left( \begin{matrix} \frac{2 n^2}{-\lambda + \sqrt{\lambda^2 - 4
n^2}}
                        & 0 \\ 0 & \frac{2n^2}{-\lambda - \sqrt{\lambda^2 - 4 n^2}} \end{matrix}
                        \right)  \\ & \times
\left( \begin{matrix} \frac{-1}{\sqrt{\lambda^2 - 4 n^2}}
                        & \frac{1}{2} - \frac{\lambda}{2 \sqrt{\lambda^2 - 4 n^2}} \\
                        \frac{1}{\sqrt{\lambda^2 - 4 n^2}}
                        & \frac{1}{2} + \frac{\lambda}{2 \sqrt{\lambda^2 - 4 n^2}}
\end{matrix} \right)\end{split}\end{equation}
where the last matrix on the right hand side is the inverse of the first.

We now make a change of variables by setting
$\lambda = 2 n \cos z$
for $z \in [0, \pi].$  This change of variables gives a bijection
between $[0,\pi]$ and $[-2n,2n]$.   Since $\|A_k\|\leq 2n$, all
our eigenvalues belong to $[-2n,2n]$ and so we can use this change
of variables to compute the eigenvalues.  Then \eqref{matrix_decomposition} becomes
\begin{align*}
\begin{pmatrix} - 2n \cos z & -n^2 \\ 1 & 0 \end{pmatrix}
& = \begin{pmatrix} \frac{n}{-\cos z + i \sin z}
                        & \frac{n}{-\cos z - i \sin z} \\ 1 & 1 \end{pmatrix}\times \\ &
\begin{pmatrix} \frac{n}{-\cos z + i \sin z}
                        & 0 \\ 0 & \frac{n}{- \cos z - i \sin z} \end{pmatrix}
\begin{pmatrix} \frac{-1}{2 n i \sin z}
                        & \frac{1}{2} - \frac{\cos z}{2 i \sin z} \\
                        \frac{1}{2 n i \sin z}
                        & \frac{1}{2} + \frac{\cos z}{2 i \sin z}\end{pmatrix}.
\end{align*}
Applying this decomposition to solve the above recursion we
obtain
\begin{align*}
\left( \begin{matrix} P_k(2 n \cos z) \\ Q_k(2 n \cos z)
\end{matrix} \right) &= \left( \begin{matrix} \frac{-n}{e^{-iz}} &
\frac{-n}{e^{iz}} \\ 1 & 1 \end{matrix} \right) \left(
\begin{matrix} \left( \frac{-n}{e^{-iz}} \right)^k & 0 \\ 0 &
\left( \frac{-n}{e^{iz}} \right)^k \end{matrix} \right)\times \\ &\qquad \left(
\begin{matrix} \frac{-1}{2 n i \sin z}
                        & \frac{1}{2} - \frac{\cos z}{2 i \sin z} \\
                        \frac{1}{2 n i \sin z}
                        & \frac{1}{2} + \frac{\cos z}{2 i \sin z}
\end{matrix} \right)
\left( \begin{matrix} 1 \\ 0 \end{matrix} \right) \\
&= \left( \begin{matrix} - \frac{1}{2 n i \sin z} \left( (-n e^{i z})^{k+1} - (-n e^{- i z})^{k+1} \right) \\
                      - \frac{1}{2 n i \sin z} \left( (-n e^{i z})^k - (-n e^{- i z})^k \right)
\end{matrix} \right).
\end{align*}
Thus, we have
\begin{align*}
F_k(\lambda, 0) &= \frac{P_k(\lambda)}{Q_k(\lambda)} \\
&= -n \left(\frac{e^{i z (k+1)} - e^{-i z (k+1)}}{e^{i z k} - e^{- i z k}}\right) \\
&= -n \left(\frac{\sin(z(k+1))}{\sin(z k)}\right).
\end{align*}

It now follows from Lemma \ref{lemma_phik_product} that
\begin{equation}\label{multiplicity_formula}
\begin{split}
\Phi_k(\lambda, 0) &= \Phi_k(2 n \cos z, 0)\\ &= (2n -2n \cos z)
\prod_{j=1}^k
   \left[ -n \left(\frac{\sin (z (j+1))}{\sin (z j)}\right) \right]^{(n-1) n^{k-j}} \\
   &= 2n (1-\cos z) (-n)^{(n-1) n^{\frac{(k-1)k}{2}}} \left( \frac{1}{\sin z} \right)^{(n-1) n^{k-1}} \\
&\qquad  \left( \prod_{j=2}^k (\sin(z j))^{(n-1)^2 n^{k-j}}
\right) (\sin (z (k+1)) )^{n-1}
\end{split}
\end{equation}
Thus, $\Phi_k(2 n \cos z, 0) = 0$ if and only if either $\cos z =
1$ or $\sin(z j) = 0$ for some $j \in \lbrace 2, \dots, k+1
\rbrace$. That is, if and only if either $z = 0$ or $zj = l \pi$
for some $j \in \lbrace 2, \dots, k+1 \rbrace$ and integer $l$.
That is, if and only if $z = 0$ or $z = \frac{l}{j} \pi$ with $j
\in \lbrace 2, \dots k+1 \rbrace$ and $1 \leq l \leq j$. The
corresponding values of $\lambda$ give the set of eigenvalues:
$$\left\lbrace 2n \right\rbrace \cup \left\lbrace 2 n \cos
\frac{p}{q} \pi \mid q \in \lbrace 2, \dots, k+1 \rbrace,
1 \leq p < q \right\rbrace.$$
 Notice
that $-2n$ is not an eigenvalue because the factor $\left(
\frac{1}{\sin z} \right)^{(c-1)n^{k-1}}$ compensates exactly for
the remaining factors.

Our next objective is to determine the multiplicities of the
eigenvalues. First we determine the multiplicity of $2n$ as an
eigenvalue. It is a basic result in Perron-Frobenius theory
\cite{ergodic,Ramanujan} that, for a connected $2n$-regular graph,
the multiplicity of $2n$ as an eigenvalue of the incidence matrix
is $1$. Hence, $2n$ has multiplicity $1$ as an eigenvalue of
$A_k$.

  The multiplicities of the remaining
eigenvalues can be computed from the formula
\eqref{multiplicity_formula}.
 Suppose $p$ and $q$ are such that $q \in \lbrace 2,
\dots, k+1 \rbrace$, \mbox{$1 \leq p < q$} and $(p,q)=1$; we wish
to calculate the multiplicity of $2n\cos \frac{p}{q}\pi$. If
$q=k+1$, then only the last term of \eqref{multiplicity_formula}
contributes to the multiplicity, so we have multiplicity $n-1$.
Suppose now $q\in \{2,\ldots,k\}$. Let $j \in \lbrace 2, \dots,
k+1 \rbrace$. Then we have $\sin \frac{p j}{q}\pi = 0$ if and only
if $q\mid pj$, that is, if and only if $q\mid j$. Thus, setting
$[r]$ to be the integer part of a real number $r$ and
$\chi_{\mathrm{Div}(k+1)}$ to be the characteristic function for
the set of divisors of $k+1$, we obtain that the eigenvalue $2 n
\cos \frac{p}{q} \pi$ has multiplicity:
\begin{align*}
& (n-1)^2 \sum_{i = 1}^{[\frac{k}{q}]} n^{k-qi} + (n-1)
\chi_{\mathrm{Div}(k+1)}(q)
\\&= n^k(n-1)^2\left(\frac{1 - n^{-q([\frac{k}{q}]
+1)}}{1-n^{-q}}-1\right) + (n-1)\chi_{\mathrm{Div}(k+1)}(q).
\end{align*}

Summing up, we have proven the following:
\begin{theorem}\label{finitecasemult}
The spectrum of the Markov operator $M_k$ from the random walk on
\Sch {\Ga} {P_k} {S} is: $$\left\lbrace 1 \right\rbrace \cup
\left\lbrace \cos \frac{p}{q} \pi \mid q \in \lbrace 2, \dots, k+1
  \rbrace, 1 \leq p < q \right\rbrace.$$  The eigenvalue $1$ has
multiplicity $1$.  For $1\leq p<q$ with $p$ and $q$ coprime, the
multiplicity is
\begin{equation*}
\#\left(\cos\frac{p}{q}\pi\right)=\begin{cases}
n-1 & \text{ if } q=k+1\\
\begin{aligned}n^k(n-1)^2\left(\frac{1 - n^{-q([\frac{k}{q}]
+1)}}{1-n^{-q}}-1\right) \\+
(n-1)\chi_{\mathrm{Div}(k+1)}(q)\end{aligned} & \text{else.}
\end{cases}
\end{equation*}\qed
 \end{theorem}

It is interesting to note that this result only depends on the
size $n$ of our finite group $G$ and not on the structure of $G$.
At this point, it would be easy to compute the Ihara zeta function
of $X_k$.  We instead wait to compute the zeta function for $X$.

Let $\phi$ be the Euler totient function, so that $\phi(q)$
denotes the number of positive integers less than or equal to $q$
and coprime to $q$ . Using Theorem \ref{finitecasemult} we obtain
a proof of a classical result from number theory \cite[Theorem
309]{Hardy}. A probabilistic interpretation can be given to this
result from our computation of the KNS spectral measure; see
Proposition \ref{knscalculated} in the next subsection.

\begin{corollary}\label{numbertheory}
Let $n\geq 2$ be an integer.  Then $$(n-1)^2\sum_{q=2}^{\infty}
\frac{\phi(q)}{n^q-1} = 1.$$
\end{corollary}
\begin{proof}
The operator $M_k$, being symmetric, has $n^k$ eigenvalues with
multiplicity.  Using Theorem \ref{finitecasemult} and observing that
the multiplicity of $\cos \frac{p}{q}\pi$, where $1\leq p<q$,
$(p,q)=1$, depends only on $q$, we obtain
\begin{equation*}
n^k = 1 + n^k(n-1)^2\sum_{q=2}^{k} \phi(q)\left(\frac {1 -
n^{-q([\frac{k}{q}] +1)}}{1-n^{-q}}-1\right)  + (n-1)\sum_{q\mid
k+1,q\neq 1}\phi(q).
\end{equation*}
Dividing both sides by $n^k$ gives:
\begin{equation*}
1 = \frac{1}{n^k} +
(n-1)^2\sum_{q=2}^{k}\phi(q)\left(\frac {1 - n^{-q([\frac{k}{q}]
+1)}}{1-n^{-q}}-1\right)  + \frac{(n-1)}{n^k}\sum_{q\mid k+1, q\neq 1}\phi(q).
\end{equation*}
Since $\phi(q)\leq q$, and so
$$\sum_{q\mid k+1, q\neq 1}\phi(q)\leq \frac{(k+1)(k+2)}{2}-1,$$
we see, by taking the limit as $k\to \infty$, that
\begin{align*}
1 &=  \lim_{k\to \infty}\ (n-1)^2\sum_{q=2}^{k}\phi(q)\left(\frac
{1 -
    n^{-q([\frac{k}{q}]  +1)}}{1-n^{-q}}-1\right)\\
 & =   \lim_{k\to \infty}\
    (n-1)^2\sum_{q=2}^{k}\phi(q)\left(\frac{n^q-n^{-q[\frac{k}{q}]}}{n^q-1}
    -1\right) \\\label{euler}
 &=  \lim_{k\to \infty}\
    (n-1)^2\left[\sum_{q=2}^{k}\frac{\phi(q)}{n^q-1} -
    \sum_{q=2}^{k}\phi(q)\left(\frac{n^{-q[\frac{k}{q}]}}{n^q-1}\right)\right].
\end{align*}
Thus, to establish the result, we need only to show that the last
term vanishes as $k\to\infty$. But
\begin{align*}
\sum_{q=2}^{k}\phi(q)\left(\frac{n^{-q[\frac{k}{q}]}}{n^q-1}\right)
&\leq k\sum_{q=2}^k
\left(\frac{n^{-q(\frac{k}{q}-1)}}{n^q-1}\right) = k\sum_{q=2}^k
\frac{1}{n^k-n^{k-q}}\\ & \leq
k(k-1)\frac{1}{n^{k-2}(n^2-1)}\longrightarrow
 0\ \text{as}\ k\to\infty.
\end{align*}
\end{proof}

\subsection{Calculating the KNS measure}
If $\psi$ is a probability measure defined on the Borel subsets of
$\R$, then the associated (cumulative) \emph{distribution
function} \cite{Billingsley} $F_{\psi}:\R\to [0,1]$ is given by
\begin{equation*}
F_{\psi}(x) = \psi((-\infty,x]).
\end{equation*}
One has that $F_{\psi}$ is non-decreasing, right continuous and
\begin{equation}\label{distributionfunc}
\lim _{x\to -\infty}F_{\psi}(x) = 0,\quad \lim_{x\to \infty}
F_{\psi}(x) =1.
\end{equation}

Conversely, if $F:\R\to [0,1]$ is a non-decreasing, right
continuous function satisfying \eqref{distributionfunc}, then
there is a unique probability measure $\psi_F$ on the Borel sets
of $\R$ such that $\psi_F((a,b]) = F(b)-F(a)$ \cite{Billingsley}.
The function $F$ can have at most countabley many discontinuity points
\cite{Billingsley}.  If the sum of the jumps of these points is $1$, 
then $\psi_F$ will be a
discrete measure supported at the discontinuity points of $F$ with
the weight at a discontinuity point equal to the amount of the
jump.

Suppose $\{\psi_k\}$ is a sequence of probability measures on
$\R$. Then  $\psi_k\to \psi$ weakly if and only if
$F_{\psi_n}(x)\to F_{\psi}(x)$ at each point of continuity of
$F_{\psi}$ \cite{Billingsley}.

Let $\mu_k$ be the measure associated to $M_k$ as per
\eqref{averagemeasure}.  The KNS spectral measure $\mu$ is the
weak limit of $\mu_k$.  We shall compute it by computing its
distribution function. To ease notation, we shall perform a change
of variables.  Let $f:[0,1]\to [-1,1]$ be given by $f(x) = \cos
\pi x$. Let us define the measures $\sigma_k$, $k\geq 0$, and
$\sigma$ on the Borel sets of $[0,1]$ by $$\sigma_k(B) =
\mu_k(f(B)),\quad \sigma(B) = \mu(f(B)).$$  Since $f$ is a
homeomorphism it follows that $\sigma$ and $\mu$ determine each
other and that $\sigma_k\to \sigma$ weakly.

To calculate $F_{\sigma}$, we define a one-parameter family of
Euler $\phi$-functions $\{\phi_x:\N\to \N\}_{x\in [0,1]}$, by
\begin{equation*}
\phi_x(q) = \left|\left\{p\in \N\mid (p,q)=1\ \text{and}\
\frac{p}{q}\leq x\right\}\right|.
\end{equation*}
Then, $\phi_0(q) =0$, $\phi(q_1) = \phi (q)$ and $\phi_x(q)$ is
non-decreasing as a function of $x$ for fixed $q\in \N$.

\begin{proposition}\label{knscalculated}
For all $x\in [0,1]$,
$$\lim _{k\to\infty} F_{\sigma_k}(x) = (n-1)^2\sum_{q=2}^{\infty}\frac{\phi_x(q)}{n^q-1}.$$
\end{proposition}
\begin{proof}
This proof is exactly like the proof of Corollary
\ref{numbertheory}, only in the right hand side of the various
equations, the role of $\phi(q)$ is taken by $\phi_x(q)$, while in
the left hand side the role of $1$ is taken by $F_{\sigma_k}(x)$
before taking limits. The same estimates apply since
\mbox{$\phi_x(q)\leq \phi(q)$}.
\end{proof}

Set $F=\lim F_{\sigma_k}$.  We know that $F=F_{\sigma}$ since
$\sigma_k\to \sigma$ weakly, but we prefer to verify directly that
$F$ is indeed a distribution function, thereby giving a direct
proof, independent of \cite{GZihara}, that the sequence of
measures $\sigma_k$ has a weak limit.

\begin{proposition}\label{convergence}
$F$ is a probability distribution function defining a discrete
measure $\sigma$ supported on the rational points of the interval
$(0,1)$. More precisely,
\begin{equation}\label{definemusigma}
\sigma = (n-1)^2\sum_{q=2}^{\infty}\left(\sum_{1\leq p<q,\
(p,q)=1} \frac{1}{n^q -1}\delta_{\frac{p}{q}}\right)
\end{equation}
where $\delta_{\frac{p}{q}}$ is a Dirac measure.
\end{proposition}
\begin{proof}
Since $\phi_x(q)$ is non-decreasing as a function of $x$ (for $q$
fixed), $F$ is clearly non-decreasing. By Corollary
\ref{numbertheory}, $F(1)=1$, while clearly $F(0) =0$. Now we show
right continuity.  It is immediate from Proposition
\ref{knscalculated} that if $\frac{p}{q}\in (0,1)$ is a rational
point, then
\begin{equation}\label{jumps}
\lim_{x\to \frac{p}{q}-}
\left(F\left(\frac{p}{q}\right)-F(x)\right) =
(n-1)^2\frac{1}{n^q-1}
\end{equation}
Hence, by Corollary \ref{numbertheory}, the sum of jumps at the
rational points is $1$.  It follows that $F(x)$ is continuous at
irrational points and the jump at a rational point $\frac{p}{q}$ is
$\frac{(n-1)^2}{n^q-1}$. From this,
\eqref{definemusigma} is immediate.
\end{proof}

Changing variables, observing that the set
$\{\cos\frac{p}{q}\pi\mid 1\leq p<q\}$ is dense in $[-1,1]$ and
using the freeness of the action in the case $G$ is abelian (in
conjunction with Theorem \ref{freevs1}, we obtain our main result.
\begin{theorem}\label{KNSfinal}
Let $G$ be a non-trivial finite group of order $n$.  Then the KNS spectral
measure $\mu$ for the Schreier graph of $\GG (\CC (G))$
with respect to a parabolic subgroup
$P$ and generators \eqref{generatorsS} is a discrete measure given by
\begin{equation}\label{spectral measure}
\mu = (n-1)^2\sum_{q=2}^{\infty}\left(\sum_{1\leq p<q, (p,q)=1}
\frac{1}{n^q -1}\delta_{\cos\frac{p}{q}\pi}\right)
\end{equation}
The following equalities of spectra of Hecke operators hold:
$$[-1,1]=\Sp {H_{\pi}}=\Sp {H_{\lambda_{\GG(\CC (G))/P}}} = \Sp
{H_{\lambda _{\GG(\CC (G))}}},$$ so the Markov operator for the
simple random walk on the Cayley graph of $\GG (\CC (G))$ has
spectrum $[-1,1]$.

In the case $G$ is an abelian group, $\GG (\CC(G)) = G\rwr \Z$ and
\eqref{spectral measure} gives the Kesten spectral measure of the
Markov operator for the simple random walk on the Cayley graph of
$G\rwr \Z$ with respect to the automaton generators.
\end{theorem}
For the case $G$ is abelian, the results for the Markov operator
were obtained in a different way by Dicks and Schick \cite{Dicks}.
Their result concerned random walks on wreath products $G\rwr \Z$ with
$G$ a finite non-trivial group.  If $G=\{g_1,\ldots,g_n\}$ and $\Z =\langle
t\rangle$, then Dicks and Schick used the symmetric generating set
\begin{equation}\label{dicksandschickgen}
S=\{tg_1,\ldots tg_n, g_1t\inv,\ldots, g_nt\inv\}.
\end{equation}
  It is easy to
check that in the case $G$ is abelian, these are the generators
obtained using the Cayley machine for $G$.  It is also straightforward
to verify that if $G$ and $H$ are finite groups of the same order $n$,
then the walks on $G\rwr \Z$ and $H\rwr \Z$ with the above generators
give rise to isomorphic Markov chains and hence have the same spectral
measure.  One can think of the system as being  
the Cayley graph of $\Z$, with at each point a lamp which can
either be off or illuminated in any one of $n-1$ different colours.  There is a
lamplighter who at each move, with equal probability, either moves to the
right and changes the lamp at his new position to off or to any of the
$n-1$ colours, or he can change the lamp where he currently is to off
or to any of the $n-1$ colours and then move to the left.  The starting
configuration has all lamps off and the lamplighter at the origin.  Hence
Theorem \ref{KNSfinal} gives the following result.

\begin{theorem}\label{KNSfinal2}
Let $G$ be a non-trivial finite group of order $n$.  Then the spectral
measure for the random walk on $G\rwr \Z$ with respect to the
generators \eqref{dicksandschickgen}
\begin{equation}\label{spectral measure2}
\mu = (n-1)^2\sum_{q=2}^{\infty}\left(\sum_{1\leq p<q, (p,q)=1}
\frac{1}{n^q -1}\delta_{\cos\frac{p}{q}\pi}\right)
\end{equation}
\end{theorem}

We can now calculate the zeta function using Theorem
\ref{KNSfinal} and \eqref{zetaintegralform}.

\begin{corollary}\label{zetacalc}
Let $G$ be a non-trivial finite group and $P$ be a parabolic
subgroup of $\GG (\CC (G))$.  Let $X = \Sch {\GG (\CC (G))} P S$,
$S$ as per \eqref{generatorsS}.  Then
\begin{align*}
\zeta_X(t) &= (1-t^2)^{-(n-1)}\times\\ &\qquad
\prod_{q=2}^{\infty}\left[\prod_{1\leq p
<q,\ (p,q)=1} \left(1 - 2nt\cos\frac{p}{q}\pi +
(2n-1)t^2)^{-(n-1)^2\frac{1}{n^q-1}}\right)\right].\end{align*}
  This
product converges for $|t|<\frac{1}{2n-1}$.  Moreover, if $G$ is
abelian, then $\zeta_X = \zeta_{G \rwr \Z}$.\qed
\end{corollary}

\section{The structure of Cayley machines of non-abelian groups}

In \cite{reset}, the second and third authors showed that for
finite \textit{abelian} groups $G$, the automata group
$\GG(\CC(G))$ is isomorphic to the wreath product $G\rwr \Z$. In
this section, we consider the case in which $G$ is not abelian,
showing that in most cases, the group $\GG(\CC(G))$ cannot be
expressed as a wreath product of any finite group with any
torsion-free group.  The following simple proposition allows us to
consider separately two different cases.

\begin{proposition}
Let $G$ be a finite group. Then either:
\begin{itemize}
\item[(i)] $G$ has a non-central element of odd order; or
\item[(ii)] $G$ is the direct product of a $2$-group and an
abelian group.
\end{itemize}
\end{proposition}
\begin{proof}

Suppose (i) does not hold, that is, that all odd order elements of
$G$ are central. We claim first that $G / Z(G)$ is a $2$-group.
Indeed, suppose not. Then some coset $g Z(G) \in G / Z(G)$ has
order an odd prime $q$. It follows that $g^n \in Z(G)$ if and only
if $q$ divides $n$, and in particular that $q$ divides the order
of $g$. Suppose the order of $g$ is $q^i r$ where $q\nmid r$. Then
$g^r$ has order $q^i$ but is not contained in $Z(G)$, which
contradicts the supposition that all odd order elements of $G$ are
central. Hence, $G / Z(G)$ is a $2$-group.

In particular, $G$ is a central extension of a nilpotent group,
and so is nilpotent. Hence, $G$ is a direct product of its Sylow
subgroups.  But for odd primes $p$, the $p$-Sylow subgroups are
central by assumption and so in particular must be abelian. Thus
$G$ is a direct product of a 2-group and an abelian group.
\end{proof}

We shall show that in the first case, our group $\GG (\CC (G))$
cannot be a wreath product of a finite group with a torsion-free
group. The second case is slightly more involved and we can only
handle the case where the $2$-group component is not nilpotent of
class $2$. In this case we again show that $\GG (\CC (G))$ cannot
embed in a wreath product of a finite group with a torsion-free group.

Fix now a finite group $G$. We consider the free monoid $G^*$ over
the elements of $G$ and write elements as bracketed,
comma-separated sequences, to avoid confusion with the
multiplication in $G$. Set $x=\CC(G)_1\inv \in \Ga$. Then we
recall from \cite{reset} that
\begin{equation}\label{resetformula}
\begin{split}
x(g_0,g_1,\ldots,g_n) &= (g_0,g_0\inv g_1,g_1\inv
g_2,\ldots,g_{n-1}\inv g_n)\\
x\inv (g_0,g_1,\ldots,g_n) &= (g_0,g_0g_1,\ldots,g_0g_1\cdots
g_n)\\
x\CC(G)_g(g_0,g_1,\ldots,g_n)&= (gg_0,g_1,\ldots, g_n)
\end{split}
\end{equation}
for every $g \in G$. It follows that the elements of the form
$x\CC(G)_g$ form a subgroup of $\GG (\CC(G))$ isomorphic to $G$.
For notational convenience, we identify $g$ with $x\CC(G)_g$ and
view $G$ as embedded in $\GG(\CC(G))$.  It is shown in
\cite{reset} that the element $x^ngx^{-n}$ has depth exactly
$n+1$.

There is an infinite sequence of words that will play an important
technical role in what follows.  Let $X=\{t_0,t_1,t_2,\ldots\}$ be
a countably infinite set, which we view as a set of variables.  In
what follows, we take the viewpoint that $X^*$ consists of terms.
If $w\in X^*$, We sometimes write $w(t_0,\ldots,t_n)$ if $w\in
\{t_0,\ldots,t_n\}^*$.  If $m_0,\ldots,m_n$ are elements of a
monoid $M$, we write $w(m_0,\ldots,m_n)$ to denote the element of
$M$ obtained by substituting $m_i$ for $t_i$.  Our sequence
$\{w_n\}$ is defined recursively as follows:
\begin{itemize}
\item  $w_0(t_0) = t_0$ and, for $n\geq 0$, \item
$w_{n+1}(t_0,t_1,\ldots,t_n, t_{n+1}) =
w_n(t_0,t_0t_1,\ldots,t_0t_1\cdots t_n)t_0t_1\cdots t_{n+1}$
\end{itemize}
Notice that $w_n$ has content $\{t_0,\ldots,t_n\}$. The first four
terms of $\{w_n\}$ are: $t_0$, $t_0t_0t_1$,
$t_0t_0t_0t_1t_0t_1t_2$ and
$t_0t_0t_0t_0t_1t_0t_0t_1t_0t_1t_2t_0t_1t_2t_3$.  Sometimes it will
be convenient to set $w_{-1}$ to be the empty word $\varepsilon$,
as the recursion still applies if we follow the usual conventions
regarding empty variable sets.

If $w\in X^*$, we denote by $|w|_{t_i}$ the number of occurrences
of $t_i$ in $w$.

\begin{lemma}\label{counting}
For $0\leq i\leq n$, $|w_n|_{t_i} = 2^{n-i}$.
\end{lemma}
\begin{proof}
The proof proceeds by induction on $n$.  For $n=0$, $w_0=t_0$ and
so the lemma holds for this case.  Suppose that the lemma holds
for $n\geq 0$.  Then $$w_{n+1}(t_0,\ldots,t_{n+1}) =
w_n(t_0,t_0t_1,\ldots,t_0\cdots t_n)t_0\cdots t_{n+1}.$$  First we
consider $1\leq i\leq n$.  In $w_n(t_0,t_0t_1,\ldots,t_0\cdots
t_n)$, there is one occurrence of $t_i$ for each occurrence of
$t_j$, $i\leq j\leq n$, in $w_n(t_0,\ldots,t_n)$.  So by
induction, we obtain $$|w_{n+1}|_{t_i} = 1+\sum_{j=i}^n
2^{n-j} = 2^{n+1-i}.$$  Clearly $|w_{n+1}|_{t_{n+1}} =1 = 2^0$.
This completes the induction, thereby establishing the lemma.
\end{proof}

The next lemma connects the sequence $\{w_n\}$ to our automata
groups.

\begin{lemma}\label{lastentry}
Let $g\in G$.  Then the last entry of $x^ngx^{-n}(g_0,\ldots,g_n)$
is $(g^{(-1)^n})^{w_{n-1}(g_0,\ldots,g_{n-1})}g_n.$
\end{lemma}
\begin{proof}
The proof is by induction on $n$.  For $n=0$, $g(g_0) = (gg_0)$,
while $(g^{(-1)^0})^{\varepsilon}g_0 = gg_0$ (recall that
$w_{-1}=\varepsilon$). Let us assume, by way of induction, that the
lemma holds for $n\geq 0$. Then
\begin{align}
x^{n+1}gx^{-(n+1)}(g_0,\ldots,g_{n+1}) & =
xx^ngx^{-n}(g_0,g_0g_1,\ldots,g_0\cdots g_{n+1})
\\ \label{secondequation} & = x(u_0,\dots,u_n,g_0\cdots g_{n+1})
\end{align}
for certain $u_i\in G$ (since $x^ngx^{-n}$ has depth $n+1$).
Moreover, we know by induction that $$u_n 
= (g^{(-1)^n})^{w_{n-1}(g_0,g_0g_1,\ldots,g_0\cdots
g_{n-1})}g_0\cdots g_n.$$   Hence, the last entry in
\eqref{secondequation} is
\begin{align*}
u_n\inv g_0\cdots g_{n+1} &= (g_0\cdots g_n)\inv
(g^{(-1)^{n+1}})^{w_{n-1}(g_0,g_0g_1,\ldots,g_0\cdots
g_{n-1})}(g_0\cdots g_n)g_{n+1}\\
&= (g^{(-1)^{n+1}})^{w_{n-1}(g_0,g_0g_1,\ldots,g_0\cdots
g_{n-1})g_0\cdots g_n}g_{n+1}\\
&=(g^{(-1)^{n+1}})^{w_{n}(g_0,\ldots,g_n)}g_{n+1},
\end{align*}
as required.
\end{proof}

Our key obstruction to embedding in wreath products is presented
by the following observation.
\begin{lemma}\label{lemma_wreath_finite_conjugacy}
Let $A=G \rwr H$ with $G$ a finite group and $H$ a torsion-free
group. Then the set of torsion elements of $A$ is the subgroup
$N=\oplus_{H} G$.  Hence every conjugacy class of $N$ is finite.
\end{lemma}
\begin{proof}
Since $A=(\oplus_{H}G)\rtimes H$ and $H$ is torsion-free, the
torsion elements of $A$ are exactly the elements of the subgroup
$N=\oplus_{H} G$. Since in the group $\oplus_{H}G$, conjugate
elements have the same support and the direct sum only contains
elements of finite support, the conjugacy classes of $N$ are
finite.
\end{proof}

\begin{theorem}
Let $G$ be a finite group with a non-central element of odd order.
Then $\GG (\CC (G))$ does not embed in the wreath product of a
finite group with a torsion-free group.
\end{theorem}
\begin{proof}
Let $g \in G$ be a non-central element of odd order. Let $h \in G$
be an element of minimal order amongst those elements that do not
commute with $g$. Let $p$ be a prime factor of the order of $h$.
Then $h^p$ has order less than that of $h$, and so commutes with
$g$. Let $v=gh\inv$.

For each $n \in \N$, we consider the element
$$\ga_n = (x^{p^n} h\inv x^{-p^n})^{-1} v (x^{p^n} h x^{-p^n}) \in \GG (\CC (G)).$$
Each such element is a conjugate of the torsion element $v$ by
another torsion element; see \eqref{defineN0}. Our objective is to
show that the $\ga_n$ are all distinct. By Lemma
\ref{lemma_wreath_finite_conjugacy}, this cannot happen in a
wreath product of a finite group with a torsion-free group, so it
will follow that $\GG (\CC (G))$ cannot embed in such a wreath
product.

To this end, we consider the action of $\ga_n$ on the word
$(1,1,\dots, 1) \in G^{p^n+1}$ and in particular compute the last
letter.  Our goal is to show that the action is non-trivial on the
last letter. Since $\ga_n$ has depth at most $p^{n}+1$, it will
then follow that $\ga_n$ has depth exactly $p^{n}+1$ and so the
various $\ga_n$ are all distinct. Using \eqref{resetformula} and
Lemma \ref{lastentry},
\begin{align*}
\ga_n(1,\ldots,1) &= x^{p^n}h\inv x^{-p^n}vx^{p^n}(h,1,\ldots,1)
\\ &= x^{p^n}h\inv x^{-p^n}\left(g,h^{-\binom{p^n}{1}},
h^{\binom{p^n}{2}}, \ldots, h^{(-1)^{p^n}
\binom{p^n}{p^n}}\right)\\
&= \left(\ldots,((h^{-1})^{(-1)^{p^n}})^{w}h^{(-1)^{p^n}}\right)
\end{align*}
where $$w=w_{p^{n}-1}\left(g,h^{-\binom{p^n}{1}},
h^{\binom{p^n}{2}}, \ldots, h^{(-1)^{p^n-1}
\binom{p^n}{p^n-1}}\right)\in G.$$

 Now it is well-known that $p$ divides $\binom{p^n}{i}$ for every
prime $p$, every $n \geq 1$ and every $1 \leq i \leq p^n-1$. In
view of the fact that $h^p$ commutes with $g$, we see using Lemma
\ref{counting} that
$$w=g^{2^{p^n-1}}h^s,\ \text{where}\ s=\sum_{i=1}^{p^n-1}
(-1)^i\binom{p^n}{i}2^{p^n-1-i}.$$ Now
\begin{equation}\label{notwreatheq1}
\begin{split}
&\qquad \quad ((h^{-1})^{(-1)^{p^n}})^{w}h^{(-1)^{p^n}} = 1\\
&\iff h^{-s}g^{-2^{p^n-1}}hg^{2^{p^n-1}}h^s  = h\\ &\iff
g^{-2^{p^n-1}}hg^{2^{p^n-1}} = h\\ &\iff hg^{2^{p^n-1}} =
g^{2^{p^n-1}}h
\end{split}
\end{equation}
But $g$ has odd order, so there is an integer $k$ such that
$(g^{2^{p^n-1}})^k=g$.  Since $g$ does not commute with $h$, we
conclude $g^{2^{p^n-1}}$ does not commute with $h$. Therefore,
none of the equalities in \eqref{notwreatheq1} hold and so $\ga_n$
acts non-trivially on the last letter of $(1,\ldots,1)\in
G^{p^n+1}$, finishing the proof.
\end{proof}

We now consider the case in which $G$ contains a $2$-group, which
is not nilpotent of class $2$.

\begin{theorem}
Let $G$ be a finite group containing a $2$-subgroup that is not
nilpotent of class $2$. Then $\GG (\CC (G))$ does not embed in a
wreath product of a finite group with a torsion-free group.
\end{theorem}
\begin{proof}
As in the previous proof, it will suffice to show that some
torsion element in $\GG(\CC(G))$ has infinitely many distinct
conjugates by other torsion elements.

Let $K$ be a $2$-group in $G$ that is not nilpotent of class $2$.
Let $g, f, h \in K$ be such that $h^{-1} f h f^{-1}$ does not
commute with $g$. Set
$$\ga_n = x^n g x^{-n} h x^n g^{-1} x^{-n}.$$  We claim that it
suffices to show that the last entry of
$$x^ngx^{-n}(1,\ldots,1,f,1,1)$$ differs from the last entry of
$$x^ngx^{-n}(h, \dots, 1, f, 1,
1)$$ for infinitely many positive integers $n$.  (Here each of the
strings we are acting on has length $n+1$). Indeed, for each $n$
such that this is true, we see that
\begin{align*}
 \ga_n[x^n g x^{-n} (1, 1, \dots, 1, f, 1, 1)] &= x^n g x^{-n} h
  (1, 1, \dots, 1, f, 1, 1) \\
&= x^n g x^{-n} (h, 1, \dots, 1, f, 1, 1)
\end{align*}
differs from
$$x^n g x^{-n} (1, 1, \dots, 1, f, 1, 1)$$
in position $n+1$. Thus, the element $\ga_n$ acts non-trivially on
the $n+1$ level; we know that $\ga_n$ can have depth at most
$n+1$, so it must have depth exactly $n+1$. So if we have
infinitely many $n$s for which this is the case, we have
conjugates of a given torsion element with arbitrarily large
depth, so there must be infinitely many of them.

By Lemma \ref{lastentry}, we have that the last entry of
$$x^ngx^{-n}(1, \dots, 1, f, 1, 1)\ \text{is}\ (g^{(-1)^n})^{w_{n-1}(1, \dots, 1, f, 1)}$$
and the last entry of
$$x^ngx^{-n}(h,1, \dots, 1, f, 1, 1)\ \text{is}\ (g^{(-1)^n})^{w_{n-1}(h, \dots, 1, f, 1)}.$$
 By Lemma \ref{counting} $|w_{n-1}|_{t_{n-2}} = 2$.   Hence
$$w_{n-1}(1, \dots, 1, f, 1) = f^2.$$ Let $k=w_{n-1}(h, \dots, 1,
f, 1)\in G$.  Then it will suffice to show that for infinitely
many $n$,
$$f^{-2} g f^{2} \neq k^{-1} gk,$$
that is, that for infinitely many $n$, $g$ does not commute with
$k f^{-2}$.

Now from the definition of $w_{n-1}$, we have
\begin{align*}
w_{n-1}(t_0, 1, \ldots,1, t_{n-2},1) &= w_{n-2}(t_0, t_0,
\ldots,t_0,t_0t_{n-2})t_0t_{n-2} \\
                               &= w_{n-3}(t_0, {t_0}^2, \dots,
                               t_0^{n-2}) t_0^{n-1} t_{n-2} t_0
                               t_{n-2}.\\
\end{align*}
Observing that $w_{n-3}(t_0, {t_0}^2, \dots, t_0^{n-2})$ must be a
power of $t_0$ and recalling that $|w_{n-1}|_{t_0}=2^{n-1}$ by
Lemma \ref{counting}, we obtain
\begin{align*}
w_{n-1}(t_0, 1, \ldots,1, t_{n-2},1) &=
t_0^{2^{n-1}-1}t_{n-2}t_0t_{n-2}
\end{align*}
Substituting $h$ for $t_0$ and $f$ for $t_{n-2}$, we obtain
$k=h^{2^{n-1}-1}fhf$, whence
$$kf^{-2} = h^{2^{n-1} - 1} f h f^{-1}.$$ But $h$ is in the
$2$-group $K$, and so for infinitely many $n$, we will have
\mbox{$h^{2^{n-1}} = 1$}. But then $k f^{-2} = h^{-1} f h f^{-1}$,
which by assumption does not commute with $g$, as required.
\end{proof}

The question of whether $\GG (\CC (G))$ is a wreath product of a
finite group and a torsion-free group remains open for a small
class of groups, namely for those groups $G$ that are direct
products of an abelian group with a $2$-group of nilpotency class
$2$. Examples include $D_4$ and the $8$-element quaternion group.

Another interesting question is that of whether, for non-abelian
$G$, the group $\GG (\CC (G))$ has bounded torsion. We conjecture
that this is never the case.

If, as we suspect, the group $\GG (\CC (G))$ is not isomorphic to
$G\rwr \Z$ for any non-abelian $G$, then one might ask whether
this wreath product arises as an automata group at all and, if so,
what conditions can be placed upon the automaton.  In particular,
can a reset automaton be found?

\bibliographystyle{amsplain}

\end{document}